\documentclass[a4paper,11pt]{article}

\usepackage{ucs}
\usepackage[utf8x]{inputenc}
\usepackage{graphicx}
\usepackage{amsfonts}
\usepackage{amssymb}
\usepackage{amsmath}
\usepackage{amsthm}
\usepackage{dsfont}
\usepackage{enumerate}
\usepackage{fullpage}
\usepackage{ifthen}
\usepackage{epic}
\usepackage{authblk}
\usepackage{textcomp}
\usepackage[all]{xy}
\usepackage{mathrsfs}
\usepackage{mathabx}
\usepackage{subcaption}
\usepackage{titlesec}

\usepackage{etoolbox}
\makeatletter
\patchcmd{\ttlh@hang}{\parindent\z@}{\parindent\z@\leavevmode}{}{}
\patchcmd{\ttlh@hang}{\noindent}{}{}{}
\makeatother

\usepackage[hypertexnames=false,colorlinks=true,linkcolor=blue,citecolor=blue]{hyperref}
\usepackage[numbers,comma,square,sort&compress]{natbib}
\usepackage[letterpaper,text={7in,9in},centering]{geometry}

\usepackage{color}
\usepackage{titlesec}
\graphicspath{{eps/}{pdf/}}

\numberwithin{equation}{section}

\newtheorem{lem}{Lemma}[section]
\newtheorem{thm}{Theorem}
\newtheorem{prop}[lem]{Proposition}
\newtheorem{cor}[lem]{Corollary}
\newtheorem{rmk}[lem]{Remark}
\newtheorem{defi}[lem]{Definition}

  {\begin{trivlist}\item[]{\bf Hypothesis #1 }\em}{\end{trivlist}}

\newcommand{\R}{\mathbb{R}}
\newcommand{\N}{\mathbb{N}}

\newcommand{\md}{\mathrm{d}}
\newcommand{\K}{\mathcal{K}}
\newcommand{\ve}{\varepsilon}
\newcommand{\de}{\delta}
\newcommand{\mz}{\mathbf{z}}
\newcommand{\mx}{\mathbf{x}}

\newcommand{\my}{\mathbf{y}}
\newcommand{\bfu}{\mathbf{u}}

\newcommand{\mZ}{\mathcal{Z}}
\newcommand{\ds}{\displaystyle}

\newcommand{\vertiii}[1]{{\left\vert\kern-0.25ex\left\vert\kern-0.25ex\left\vert
          #1
        \right\vert\kern-0.25ex\right\vert\kern-0.25ex\right\vert}}

\makeatletter\@addtoreset{figure}{section}\makeatother

\makeatletter \@addtoreset{equation}{section} \makeatother

\newcommand{\Supp}{\mathop{\mathrm{Supp}}}

\definecolor{darkblue}{rgb}{0,0,1}

\newenvironment{Proof}%
{\begin{trivlist} \item[]{\bf Proof. }}%
{\hspace*{\fill}$\rule{.4\baselineskip}{.4\baselineskip}$\end{trivlist}}

\title{\textcolor{black}{Asymptotic} limit of a spatially-extended \textcolor{black}{mean-field} FitzHugh-Nagumo model}
\author[1]{Joachim Crevat\footnote{ \texttt{joachim.crevat@math.univ-toulouse.fr}}}
\affil[1]{Institut de Math\'ematiques de Toulouse ; UMR5219, Universit\'e de Toulouse ; UPS IMT, F-31062 Toulouse Cedex 9 France}
\begin{document}
\maketitle

\begin{abstract}
We consider a spatially extended \textcolor{black}{mean-field} model of a FitzHugh-Nagumo neural network, with a rescaled interaction kernel. Our main purpose is to prove that its \textcolor{black}{asymptotic limit} in the regime of strong local interactions \textcolor{black}{converges toward a system of reaction-diffusion equations} taking account for the average quantities of the network. Our approach is based on a \textcolor{black}{modulated energy} argument, to compare the macroscopic quantities computed from the solution of the \textcolor{black}{transport} equation, and the solution of the limit system. The main difficulty, compared to the literature, lies in the need of regularity in space of the solutions of the limit system and a careful control of an internal nonlocal \textcolor{black}{dissipation}.
\end{abstract}

\begin{center}
\begin{tabular}{ll}
&
\textbf{\textit{Keywords: }} \textcolor{black}{Asymptotic} limit, \textcolor{black}{modulated energy}, FitzHugh-Nagumo, neural network. 
\\ \,&\\
&
\textbf{\textit{Mathematics Subject Classification (2010): }} 35Q92, 35K57, 82C22, 92B20
\end{tabular}
\end{center}

\section{Introduction}

{\color{black}



In this article, we consider the following nonlocal transport equation for all $\ve>0$, $t>0$, $\mx\in\R^d$ with $d\in\{1,2,3\}$, and $\bfu=(v,w)\in\R^2$:
\begin{equation}\label{eq:f}
\left\{\begin{array}{l}
\partial_t f^\ve(t,\mx,\bfu) \,+\,\partial_v\left[ f^\ve(t,\mx,\bfu)\,\left( N(v) - w - \K_\ve[f^\ve](t,\mx,v) \right) \right] \,+\,\partial_w\left[ f^\ve(t,\mx,\bfu)\,A(v,w) \right] \,=\,0,\\ \,  \\
\K_\ve[f^\ve](t,\mx,v) \,:=\,\dfrac{1}{\ve^{d+2}}\ds\iint\Psi\left(\dfrac{\|\mx-\mx'\|}{\ve}\right)\,(v-v')\,f^\ve(t,\mx',\bfu')\,\md\mx'\,\md \bfu' ,\\ \, \\
A(v,w)\,:=\,\tau\,(v-\gamma\,w), \\ \, \\
f^\ve|_{t=0} \,=\,f^\ve_0,
\end{array}\right.
\end{equation}
using the shorthand notation $\bfu'=(v',w')$. It was derived in \cite{CREmf} as the mean-field limit of a spatially-organized FitzHugh-Nagumo (FHN) system for the modeling of a finite size neural network, as the number of neurons goes to infinity (see \cite{FIT,NAG} for the single neuron model). Thus, in this framework, each neuron is characterized by three quantities: $v\in\R$ its membrane potential, $w\in\R$ which is called adaptation variable, and $\mx\in\R^d$ its spatial position in the network. Furthermore, $f^\ve(t,\mx,\bfu)$ is the density function of finding neurons at time $t$ at position $\mx$ and at the electrical state $\bfu=(v,w)$ within the cortex.


Let us clarify the notations. First of all $\tau$ and $\gamma$ are non-negative constant. Here, we consider small values of $\tau>0$ to account for the slow evolution of the adaptation variable. Then, the purpose of the function $N$ is to model the excitability of each neuron. In the rest of this article, as in \cite{CAR,CS,HUP10,HUP13,JON,QUI}, we consider the following cubic nonlinearity 
\begin{equation}\label{N}
N(v) \,:=\, v\,(1-v)\,(v-a),\quad v\in\R,
\end{equation}
where $a\in (0,1)$ is fixed. Such a nonlinearity satisfies the following properties for four fixed positive constants $\kappa_1$, $\kappa_1'$, $\kappa_2$ and $\kappa_3$:
\begin{equation}\label{hyp:N}
\left\{\begin{array}{ll}
v\,N(v)\,\leq\,\kappa_1\,|v|^2 \,-\, \kappa_1'\,|v|^4\, , & v\in\R, \\ \, \\
(v-u)\,(N(v)-N(u)) \, \leq\kappa_2\,|v-u|^2 \, , & v,u\in\R, \\ \, \\
\left| N(v)-N(u)  \right| \, \leq \, \kappa_3\,|v-u|\,\left( 1+|v|^2+|u|^2 \right) \, , & v,u\in\R.
\end{array}\right.
\end{equation}
In the rest of this article, we consider a small scaling parameter $\ve>0$ and
\begin{equation}\label{eq:phi}
\Psi_\ve(\|\mx\|)\,:=\, \dfrac{1}{\ve^d}\,\Psi\left(\dfrac{\|\mx\|}{\ve}\right),\quad\mx\in\R^d,\quad \ve>0,
\end{equation}
where $\Psi:\R\rightarrow\R^+$ is a connectivity kernel modeling the influence of the distance between two neurons on their interactions. Notice that in this model, interactions are only modulated by distance and the electrical voltage between neurons. In the following, we assume that $\Psi$ satisfies the following assumptions:
\begin{equation}\label{hyp:Psi}
\Psi> 0 \text{ a.e.},\quad \ds\int\Psi(\|\my\|)\,\md\my \,=\, 1,\quad 0\,<\, \sigma\,=\,\ds\int\Psi(\|\my\|)\,\dfrac{\|\my\|^2}{2}\,\md\my\,<\,\infty .
\end{equation}
Let us discuss the hypotheses on $\Psi$. First, the choice of a positive and integrable connectivity kernel implies that we only consider activatory interactions, and that neurons which are far away from each other have few interactions. Moreover, we restrict our framework to positive functions only for technical reasons. Then, the assumptions of finite moments and symmetry are crucial hypotheses to derive the asymptotic limit as $\ve$ goes to $0$. Indeed, the finite moment assumption determines the type of spatial diffusion we get in the limit, and as we shall see later, the symmetry assumption yields that all the moments of $\Psi$ of odd order are $0$. A typical example of admissible connectivity kernel in our framework is a Gaussian function. 

The main purpose of this article is to derive a macroscopic model of the neural network from the transport equation \eqref{eq:f}, which accounts for the evolution of the average membrane potential of neurons at each position $\mx$. Yet, the macroscopic values computed from $f^\ve$ a solution to \eqref{eq:f} do not satisfy a closed system of equations. Indeed, let us consider $\rho^\ve$ the average density of neurons, $V^\ve$ the average membrane potential and $W^\ve$ the average adaptation variable, defined through
\begin{equation}\label{def:macro}
\left\{\begin{array}{l}
\rho^\ve(t,\mx) \,:=\, \ds\int_{\R^2} f^\ve(t,\mx,v,w)\,\md v\,\md w, \\ \, \\
\rho^\ve(t,\mx)\,V^\ve(t,\mx) \,:=\, \ds\int_{\R^2} v\,f^\ve(t,\mx,v,w)\,\md v\,\md w, \\ \, \\
\rho^\ve(t,\mx)\,W^\ve(t,\mx) \,:=\, \ds\int_{\R^2} w\,f^\ve(t,\mx,v,w)\,\md v\,\md w.
\end{array}\right.
\end{equation}
We directly get from \eqref{eq:f} the conservation of mass, that is for all $t>0$, $\rho^\ve(t,\cdot) = \rho^\ve(0,\cdot) = \rho_0^\ve$. Then, these macroscopic quantities formally satisfy the following system for all $t>0$ and $\mx\in\R^d$:
\begin{equation}\label{eq:unclosed}
\left\{\begin{array}{l}
\begin{aligned}\partial_t(\rho_0^\ve\,V^\ve)(t,\mx) \,&-\, \dfrac{1}{\ve^2}\ds\int_{\R^d}\Psi_\ve(\|\mx-\mx'\|)\,\left( V^\ve(t,\mx')-V^\ve(t,\mx) \right)\,\rho_0^\ve(\mx)\,\rho_0^\ve(\mx')\,\md\mx' \\&=\, \ds\int_{\R^2}N(v)\,f^\ve(t,\mx,v,w)\,\md v\,\md w \,-\, \rho_0^\ve\,W^\ve(t,\mx) ,\end{aligned} \\ \, \\
\partial_t(\rho_0^\ve\,W^\ve)(t,\mx) \,=\, \rho_0^\ve(\mx)\,A\left( V^\ve(t,\mx) , W^\ve(t,\mx) \right).
\end{array}\right.
\end{equation}
The problem comes from the fact that the function $N$ is not linear. The whole interest of the scaling of the nonlocal term in \eqref{eq:f} is to consider the regime of strong local interactions, in order to circumvent this issue. Indeed, in the asymptotic limit $\ve\rightarrow 0$, the interaction kernel $\Psi_\ve$ converges towards a Dirac distribution. From a biological viewpoint, this seems definitely justified, since at the macroscopic scale, if we consider a whole portion of the cortex, any neuron seems to interact only with its very closest neighbours, so the nonlocal effect is insignificant. 


In the spirit of \cite{CFF}, we expect the solution of the transport equation \eqref{eq:f} to converge towards a Dirac distribution in $v$, centered in a function $V(t,\mx)$ solution of a reaction-diffusion equation, which will therefore provide a macroscopic description of the FHN neural network.
}


Now, let us formally derive the behavior of a solution $f^\ve$ of the \textcolor{black}{transport} equation \eqref{eq:f} in the limit $\ve\rightarrow 0$. Using the change of variable $\my=(\mx-\mx')/\ve$, we formally have for all $t\in[0,T]$, $\mx\in\R^d$ and $\bfu=(v,w)\in\R^2$:
$$ f^\ve(t,\mx,\bfu)\,\K_\ve[f^\ve](t,\mx,v) \,=\, \ve^{-2}\ds\iint\Psi(\|\my\|)\,(v-v')\,f^\ve(t,\mx-\ve\,\my,\bfu')\,f^\ve(t,\mx,\bfu)\,\md\my\,\md \bfu'.$$
Moreover, using a Taylor expansion, we get that for all $t\in[0,T]$, $\bfu'\in\R^2$ and $(\mx,\my)\in\R^{2d}$,
\begin{equation}\label{eq:taylor}
f^\ve(t,\mx-\ve\my,\bfu')\,=\,f^\ve(t,\mx,\bfu') \,-\,\ve\,\nabla_\mx f^\ve(t,\mx,\bfu')\cdot\my\,+\,\dfrac{\ve^2}{2}\,\my^T\cdot \nabla_\mx^2 f^\ve(t,\mx,\bfu')\cdot\my \,+\,... \,,
\end{equation}
\textcolor{black}{where $\nabla_\mx^2 f^\ve$ is the Hessian matrix of $f^\ve$ with respect to $\mx$. }We resume the computation of $f^\ve\,\K_\ve[f^\ve]$ using the expansion \eqref{eq:taylor}. \textcolor{black}{According to the assumptions \eqref{hyp:Psi} satisfied by $\Psi$, since the connectivity kernel $\Psi$ is radially symmetric, we can simplify all the terms of odd order, which gives:}
\begin{eqnarray}
f^\ve(t,\mx,\bfu) \,\K_\ve[f^\ve](t,\mx,v) &=& \ve^{-2}\ds\int (v-v')\,f^\ve(t,\mx,\bfu')\,f^\ve(t,\mx,\bfu)\,\md \bfu'   \nonumber\\
&&\,+\,\sigma \,\ds\int (v-v')\,\Delta_\mx f^\ve(t,\mx,\bfu')\,f^\ve(t,\mx,\bfu)\,\md \bfu' \,+\, R_\ve(t,\mx,\bfu),  \label{eq:fKf}
\end{eqnarray}
where $R_\ve(t,\mx,\bfu)$ gathers all the remaining terms. Assume that the solution $f^\ve$ converges towards a distribution $f$ in some weak sense, and that $R_\ve$ tends towards $0$ as $\ve$ goes to $0$. In the following, we define the \textcolor{black}{triple} $(\textcolor{black}{\rho_0},V,W)$ of macroscopic quantities associated to $f$ similarly as in \eqref{def:macro}\textcolor{black}{, where $\rho_0$ does not depend on time}. Then, we insert the expansion \eqref{eq:fKf} into the \textcolor{black}{transport} equation \eqref{eq:f}, and we identify the orders $-2$ and $0$ of $\ve$ as $\ve$ tends to $0$. \textcolor{black}{This} formally leads to the following equations satisfied by $f$ in the sense of distributions for $t>0$, $\mx\in\R^d$, $(v,w)\in\R^2$:
\begin{eqnarray}
\mathcal{O}(\ve^{-2}): \quad &\ds\iint (v-v')\,f(t,\mx,v',w')\,f(t,\mx,v,w)\,\md v'\,\md w' \,=\,0,\label{eq:1}  \\ \, \nonumber \\
\mathcal{O}(1): \quad &\partial_t f \,+\,\partial_v\left[ f\,\left(N(v)-w\right) \right] \,+\,\partial_w\left[ f\,A(v,w) \right] \,-\,\sigma \,\partial_v\left[ f\,(\Delta_\mx\rho_0\,v-\Delta_\mx(\rho_0\,V)) \right] \,=\, 0.\label{eq:3}
\end{eqnarray}
The equation \eqref{eq:1} implies that $f$ is proportional to a Dirac mass in $v$, that is for all $t\in[0,T]$, $\mx\in\R^d$ and $(v,w)\in\R^2$, 
\begin{equation}\label{eq:mono}
f(t,\mx,v,w)\,=\,F(t,\mx,w)\,\otimes\,\delta_{V(t,\mx)}(v),
\end{equation}
where we define
$$ F(t,\mx,w)\,:=\,\ds\int_\R f(t,\mx,v,w)\,\md v, \quad \rho_0(\mx)\,W(t,\mx)\,:=\, \int_\R w\,F(t,\mx,\md w). $$
{\color{black}Indeed, \eqref{eq:1} yields that $\left(v-V(t,\mx)\right)\,\rho_0(\mx)\,f(t,\mx,v,w)\,=\,0$ in the sense of distributionssince for all $\phi\in\mathscr{C}^\infty_c(\R^{d+2})$ and $t>0$, we get: 
\begin{align*}
\left|\ds\iint f(t,\mx,v,w)\,\rho_0(\mx)\,\phi(\mx,v,w)\,\right.&\left.\md\mx\,\md v\,\md w \,-\, \ds\int F(t,\mx,w)\,\rho_0(\mx)\,\phi(\mx,V(t,\mx),w)\,\md\mx\,\md w\right| \\
&\leq\,\|\nabla_v \phi\|_{L^\infty}\ds\iint_{\text{Supp}(\phi)} f(t,\mx,v,w)\,\rho_0(\mx)\,\left|v-V(t,\mx)\right|\,\md\mx\,\,\md v\,\md w ,
\end{align*}
}
which is $0$. This concludes the justification of our claim \eqref{eq:mono}.
Therefore, we deduce from \eqref{eq:3}-\eqref{eq:mono} that the \textcolor{black}{pair} $(V,F)$ satisfies the following system for $t>0$, $\mx\in\R^d$ and $w\in\R$:
\begin{equation}\label{eq:lim}
\left\{\begin{array}{l}
\partial_t  F  \,+\,\partial_w\left( A(V,w) \,F \right) \,=\, 0, \\ \, \\ 
\partial_t \left( \rho_0 \, V \right) \,-\,\sigma \left[ \rho_0\,\Delta_\mx\left( \rho_0\,V \right) \,-\, \left(\Delta_\mx \rho_0\right) \, \rho_0\,V \right] \,=\, \rho_0\,N(V) \,-\,\rho_0\,W,   
\end{array}\right.
\end{equation}
Therefore, the limit \textcolor{black}{triple} $\mZ\,:=\,(\rho_0,\rho_0\,V,\rho_0\,W)$ is expected to satisfy the reaction-diffusion system
\begin{equation}\label{eq:rhoFHNlocal}
\partial_t \mZ \,=\, \rho_0\,\begin{pmatrix} 0 \\ N(V)-W+\sigma \big[ \Delta_\mx(\rho_0\,V) \,-\, \Delta_\mx\rho_0 \, V \big] \\ A(V,W)  \end{pmatrix}.
\end{equation} 

Let us discuss the structure of the system \eqref{eq:rhoFHNlocal}. It is worth noticing that it is not well-defined for $\mx\in\R^d$ if $\rho_0(\mx)=0$. In the case $\rho_0>0$, \textcolor{black}{it} reduces to the \textcolor{black}{FHN} reaction-diffusion system for $t>0$ and $\mx\in\R^d$:
\begin{equation}\label{eq:FHNlocal0}
\left\{\begin{array}{l}
\partial_t V(t,\mx) \,-\,\sigma \,\big[\Delta_\mx \left( \rho_0\,V \right)(t,\mx) \,-\,\Delta_\mx\rho_0(\mx)\,V(t,\mx)\big] \,=\, N(V(t,\mx)) \,-\,W(t,\mx),   \\ \, \\
\partial_t W(t,\mx) \,=\, A(V(t,\mx),W(t,\mx)).
\end{array}\right.
\end{equation}

{\color{black}
We stress that the rescaling $\ve^{-2}$ in the nonlocal term in \eqref{eq:f} is crucial. From a biological viewpoint, we only need it to go to infinity as $\ve$ goes to $0$ to provide strong interactions, but from a mathematical viewpoint, this power $-2$ is the only one which enables us to get local interactions in the limit $\varepsilon\rightarrow0$. We refer to \cite{BAT,BATES}, in which the authors use the same rescaling to prove the convergence of the traveling wave solutions of a nonlocal bistable equation towards solutions of a standard local one.
}

The system \eqref{eq:FHNlocal0} has been extensively studied, namely regarding the formation and propagation of traveling fronts and pulses for example (see e.g. \cite{CS,CAR,JON}). We also mention some recent works on the \textcolor{black}{FHN} system in the discrete case \cite{HUP10,HUP13}. The contribution of our article is thus to prove that the system \eqref{eq:lim}, which can be reduced under some assumptions to the \textcolor{black}{FHN} reaction-diffusion system \eqref{eq:FHNlocal0}, is a macroscopic description of the neural network.

\textcolor{black}{Our approach to prove the \textcolor{black}{asymptotic limit} follows ideas from \cite{CFF,FIG,KAN,KAR}, where the authors derive macroscopic equations using a relative entropy argument, as developed in the works by Dafermos \cite{Dafermos} and Di Perna \cite{Diperna} for conservation laws. Yet, the entropy used in the \cite{CFF} for instance is the standard energy functional, so the associated relative entropy actually corresponds to the notion of modulated energy as introduced in \cite{BRE} to prove the convergence of the Vlasov-Poisson system in the quasi-neutral regime towards the incompressible Euler equations. This modulated energy argument consists in estimating the modulation of the energy functional with the solution of the limit equation. It was developed precisely to derive asymptotic limits of kinetic equations without velocity transport, as in our framework, or in \cite{BRE2} for instance.
}

The specificity of our problem is the absence of noise \textcolor{black}{in the considered neural network}, which implies the absence of a Laplace operator in $v$ in the \textcolor{black}{mean-field} equation \eqref{eq:f}. Hence, without the regularizing effect of noise, the solution $f^\ve$ of \eqref{eq:f} converges towards a Dirac distribution in $v$, which prevents us from using \textcolor{black}{the decay of} a usual entropy of type $f\,\log(f)$ as in \cite{KAR}. As in \cite{CFF,FIG,KAN}, we rather focus on the evolution of moments of second order in $v$ and $w$ of $f^\ve$. Then, we encounter two main difficulties. The first one comes from the term $\partial_v(f^\ve\,N(v) )$ in \eqref{eq:f}, which introduces moments of $f^\ve$ in $v$ of higher order. Similarly as in \cite{CFF}, it will be sufficient to control moments of $f^\ve$ of fourth order to circumvent this problem. Then, unlike \cite{CFF}, we have to precisely determine the regularity in space we need for the solutions of the limit system \eqref{eq:lim} to estimate the \textcolor{black}{modulated energy}.

The problem tackled in this article lies within a rich literature of works aiming to establish a rigorous link between a mesoscopic model and a macroscopic characterization of large neural networks. For example, in \cite{PER}, the authors derived the Integrate-and-Fire model as the macroscopic limit of a voltage-conductance kinetic system. Then, the article \cite{BOS} studied the mean-field limit of a Hodgkin-Huxley neural network, and proved the existence of synchronized dynamics for the microscopic and macroscopic models. If we focus on the \textcolor{black}{FHN} model, a recent example is the work of \cite{QUI}, which investigates the large coupling limit of a kinetic description of a noisy \textcolor{black}{FHN} neural network, with uniform conductance, and highlight the emergence of clamping or synchronization between neurons in the limit. We also mention \cite{LUCON14}, in which the authors proved the mean-field limit of noisy spatially-extended \textcolor{black}{FHN} network and propagation of chaos, so that in the limit, the electrical state of each neuron can be modeled by a stochastic differential equation. This last result has been extended in \cite{LUCON18} to the mean-field limit of spatially-extended  \textcolor{black}{FHN} neural networks on random graphs, to prove the convergence towards a nonlocal reaction-diffusion system.

\paragraph{Outline of the paper.} The rest of this paper is organised as follows. In Section \ref{sec:mainresult}, we present our hypotheses, and our results on the existence of solutions to the \textcolor{black}{transport} equation \eqref{eq:f} and to the limit equation \eqref{eq:lim}, and our main result about the \textcolor{black}{asymptotic limit} from \eqref{eq:f} to \eqref{eq:lim}. Then, in Section \ref{sec:apriori}, we prove some \textit{a priori} estimates which will be key arguments for the proof of our main result. Then, Section \ref{sec:proof} is devoted to the \textcolor{black}{modulated energy} estimate and the proof of our main result. Finally, in Section \ref{sec:eqlim}, we study the well-posedness of the limit equation \eqref{eq:lim}, constructing a solution from a \textcolor{black}{pair} $(V,W)$ satisfying the reaction-diffusion system \eqref{eq:FHNlocal0} in a weak sense.

\section{Main result}\label{sec:mainresult}

In this section, we state our main result on the \textcolor{black}{asymptotic limit} of a weak solution $\left(f^\ve\right)_{\ve>0}$ of the \textcolor{black}{transport} equation \eqref{eq:f} towards a solution $(V,F)$ of the reaction-diffusion equation \eqref{eq:lim}. Before that, we have to precisely define our notion of solutions of the \textcolor{black}{transport} model \eqref{eq:f} and of the system \eqref{eq:lim}. 

\subsection{Existence of a weak solution of the \textcolor{black}{transport} equation}

In this subsection, we focus on the well-posedness of the \textcolor{black}{mesoscopic} model. First, let us specify our notion of weak solution of the \textcolor{black}{transport} equation \eqref{eq:f}.
\begin{defi}\label{defi:weaksolfeps}
We say that $f^\ve$ is a weak solution of \eqref{eq:f} with initial condition $f_0^\ve\geq 0$ if for any $T>0$, 
$$f^\ve\in\mathscr{C}^0\left([0,T],L^1(\R^{d+2})\right)\cap L^\infty\left( (0,T)\times\R^{d+2} \right),$$
and for any $\varphi\in\mathscr{C}^\infty_c([0,T)\times\R^{d+2})$, the following weak formulation of \eqref{eq:f} holds, 
\begin{equation}\label{eq:weak}
\ds\int_0^T\int f^\ve\left[ \partial_t \varphi + \left( N(v) - w - \K_\ve[f^\ve] \right)\,\partial_v\varphi + A(v,w)\,\partial_w\varphi  \right]\,\md\mz\,\md t \,+\, \ds\int f_0^\ve(\mz)\,\varphi(0,\mz)\,\md\mz \,=\,0,
\end{equation}
where $\mz = (\mx,v,w)\in\R^{d+2}$.
\end{defi}
In the rest of this paper, if $f^\ve$ is a weak solution of \eqref{eq:f}, we note $\mZ^\ve\,:=\,(\rho_0^\ve,\rho_0^\ve\,V^\ve,\rho_0^\ve\,W^\ve)$ the \textcolor{black}{triple} of macroscopic quantities computed from $f^\ve$ as in \eqref{def:macro}. Then, let us quote our result of existence and uniqueness of a weak solution to the \textcolor{black}{transport} equation \eqref{eq:f}, whose proof can be found in Proposition 2.2 from \cite{CFF}.
\begin{prop}\label{prop:eqf}
Let $\ve>0$. We consider a connectivity kernel $\Psi$ satisfying \eqref{hyp:Psi}, and an initial data $f_0^\ve$ such that
\begin{equation}\label{hyp:f0ve}
f_0^\ve\geq 0, \quad \|f_0^\ve\|_{L^1(\R^{d+2})}\,=\,1, \quad \textcolor{black}{f_0^\ve,\,\nabla_\bfu f_0^\ve} \in L^\infty(\R^{d+2}) ,
\end{equation}
where $\mathbf{u}=(v,w)$, and there exists a positive constant $R_0^\ve>0$ such that for all $\mx\in\R^d$, 
\begin{equation}\label{hyp:suppf0ve}
\Supp(f_0^\ve(\mx,\cdot))\,\subseteq\,B(0,R_0^\ve)\,\subset\,\R^2.
\end{equation}
Then, for any $T>0$, there exists a unique non-negative weak solution $f^\ve$ of \eqref{eq:f} in the sense of Definition \ref{defi:weaksolfeps}, which is compactly supported in $\mathbf{u}=(v,w)\in\R^2$.
\end{prop}

\begin{rmk}\
\begin{enumerate}
\item In the following, since the conservation of the $L^1$ norm of a weak solution $f^\ve$ of \eqref{eq:f} holds, we get that for all $t\in[0,T]$, $\|f^\ve(t,\cdot)\|_{L^1(\R^{d+2})}\,=\,\|f^\ve_0\|_{L^1(\R^{d+2})}\,=\,1.$
\item It is worth noticing that we do not need the weak solution $f^\ve$ of \eqref{eq:f} to be differentiable in space.
\end{enumerate}
\end{rmk}

\subsection{Existence of a solution of the limit system}

Our purpose is to prove the existence of a solution to the system \eqref{eq:lim}. We proceed in two steps: first, we prove the existence and uniqueness of a solution to the \textcolor{black}{FHN} reaction-diffusion system \eqref{eq:FHNlocal0}, and then we construct a solution to the system \eqref{eq:lim} from the solution to \eqref{eq:FHNlocal0}. Thus, we start by studying the following Cauchy problem for a given initial data $(V_0,W_0)$:
\begin{equation}\label{eq:FHNlocal}
\left\{\begin{array}{l}
\partial_t V(t,\mx) \,-\,\sigma \,\big[\rho_0(\mx)\,\Delta_\mx V(t,\mx) \,+\,2\,\nabla_\mx\rho_0(\mx)\cdot\nabla_\mx V(t,\mx)\big] \,=\, N(V(t,\mx)) \,-\,W(t,\mx),   \\ \, \\
\partial_t W(t,\mx) \,=\, A(V(t,\mx),W(t,\mx)), \\ \, \\
V|_{t=0} \,=\, V_0,\quad W|_{t=0} \,=\, W_0.
\end{array}\right.
\end{equation}
Before stating our existence result, we have to precisely define the notion of weak solution to the reaction-diffusion system \eqref{eq:FHNlocal}. Since all the mathematical difficulties come from the first equation in \eqref{eq:FHNlocal}, we consider the second component $W$ as a reaction term.
\begin{defi}\label{defi:weaksol}
For any $T>0$ and any given initial data $V_0$, $W_0\in H^2(\R^d)$, the \textcolor{black}{pair} $(V,W)$ is a weak solution of \eqref{eq:FHNlocal} on $[0,T]$ with initial data $(V_0,W_0)$ if $V\in L^\infty([0,T],H^2(\R^d))\cap W^{1,\infty}([0,T],L^2(\R^d))$, and $(V,W)$ verifies for all $\varphi\in H^{1}(\R^d)$, for all $t\in[0,T]$ and almost every $\mx\in\R^d$:
\begin{equation}\label{eq:weaklim}
\left\{\begin{array}{l}
\ds\int\partial_t V\,\varphi\,\md\mx \,=\,-\sigma\ds\int\rho_0\,\nabla_\mx V\cdot\nabla_\mx\varphi \,\md\mx \,+\, \sigma\ds\int \nabla_\mx\rho_0\cdot\nabla_\mx V \,\varphi\,\md\mx \,+\,\ds\int\left( N(V)-W \right)\,\varphi\,\md\mx, \\ \, \\
W(t,\mx) \,=\,  e^{-\tau\,\gamma\,t}W_0(\mx) \,+\,\tau\ds\int_0^t e^{-\tau\,\gamma\,(t-s)}V(s,\mx)\,\md s, \\ \, \\
V|_{t=0}\,=\, V_0.
\end{array}\right.
\end{equation}
\end{defi}
Now, we can give the details of our result of existence and uniqueness for the reaction-diffusion system \eqref{eq:FHNlocal}.
\begin{prop}\label{prop:eqlim}
Let $\Psi$ be a connectivity kernel satisfying \eqref{hyp:Psi}. Consider an initial data $\rho_0$ such that
\begin{equation}\label{hyp:rho0}
\rho_0\geq 0,~~\rho_0\in\mathscr{C}^{3}_b(\R^d),  ~~ \|\rho_0\|_{L^1}\,=\,1,
\end{equation}
and also consider an initial data $(V_0,W_0)$ such that
\begin{equation}\label{hyp:V0}
V_0,\,W_0\in H^2(\R^d).
\end{equation}
Then, for all $T>0$, there exists a unique function $(V,W)$ weak solution of the reaction-diffusion equation \eqref{eq:FHNlocal} on $[0,T]$ in the sense of Definition \ref{defi:weaksol} such that
$$V,\,W\in L^\infty\left([0,T],H^2(\R^d)\right)\,\cap\,\mathscr{C}^0\left([0,T],H^{1}(\R^d)\right).$$
\end{prop}

We postpone the proof of Proposition \ref{prop:eqlim} to Section \ref{sec:eqlim}. Our strategy consists in approximating the reaction-diffusion system \eqref{eq:FHNlocal} with an uniformly parabolic system, and then to pass to the limit. 

It remains to deduce from this proposition the existence of a solution to the limit system \eqref{eq:lim}. In order to work with a well-defined system on $\R^d$, we make the convention that for all $\mx\in\R^d$ such that $\rho_0(\mx)=0$, $V(\cdot,\mx)=0$ and $F(\cdot,\mx,\cdot)=0$. From a modeling viewpoint, it means that there is no electrical activity wherever there is no neuron. Hence, we are led to study the following Cauchy problem for any given initial data $(V_0,F_0)$, for all $t>0$, $\mx\in\R^d$ and $w\in\R$:
\begin{equation}\label{eq:lim2}
\left\{\begin{array}{l}
\partial_t  F  \,+\,\partial_w\left( A(V,w) \,F \right) \,=\, 0, \\ \, \\ 
\partial_t \left( \rho_0 \, V \right) \,-\,\sigma \big[ \rho_0\,\Delta_\mx\left( \rho_0\,V \right) \,-\, \left(\Delta_\mx \rho_0\right) \, \rho_0\,V \big] \,=\, \rho_0\,N(V) \,-\,\rho_0\,W,   \\ \, \\
V(t,\mx)\,=\,0, \quad t>0 \text{ and } \mx\in\R^d\backslash \text{Supp}_\text{ess}(\rho_0), \\ \, \\
W(t,\mx)\,=\,\left\{\begin{array}{l l} \dfrac{1}{\rho_0(\mx)}\ds\int_\R w\,F(t,\mx,\md w) \,&\text{if } \rho_0(\mx)>0, \\ \, \\ 0 \,&\text{else}, \end{array} \right.    \\ \, \\ 
V|_{t=0}\,=\,V_0,\quad F|_{t=0}\,=\,F_0, \quad \rho_0(\mx) \,=\,\ds\int_\R F_0(\mx,w)\,\md w.
\end{array}\right.
\end{equation}
To conclude this subsection, we state our notion of solution to the limit system \eqref{eq:lim2}, and then our result of existence and uniqueness of a solution. In the rest of this article, we denote by $\mathcal{M}(\R^{d+1})$ the set of non-negative Radon measures on $\R^{d+1}$.
\begin{defi}\label{defi:sol}
For any $T>0$, and any initial data $V_0\in H^2(\R^d)$ and $F_0\in\mathcal{M}(\R^{d+1})$ satisfying
$$\ds\int_{\R^{d+1}} |w|^2\,F_0(\md\mx,\md w) \,<\, +\infty,\quad \rho_0\,:=\,\ds\int_\R F_0(\cdot,\md w) \,\in L^1(\R^d),$$
we say that $(V,F)$ is a solution of \eqref{eq:lim2} if $F$ is a measure solution of the first equation in \eqref{eq:lim}, that is for all $\varphi\in\mathscr{C}^\infty_c(\R^{d+1})$, for all $t\in[0,T]$,
\begin{equation}
\dfrac{\md}{\md t}\ds\int_{\R^{d+1}}\varphi(\mx,w)\,F(t,\md\mx,\md w) \,-\, \ds\int_{\R^{d+1}} A(V(t,\mx),w)\,\partial_w\varphi\,F(t,\md \mx,\md w) \,=\,0,
\end{equation}
and $V\in L^\infty([0,T],H^2(\R^d))\cap W^{1,\infty}([0,T],L^2(\R^d))$ satisfies for all $\varphi\in H^1(\R^d)$ and all $t\in[0,T]$, 
\begin{equation}
\left\{\begin{array}{l}
\ds\int\partial_t(\rho_0\,V)\,\varphi\,\md\mx \,=\,\sigma\ds\int\left[(\rho_0\,V)\nabla_\mx\rho_0-\rho_0\,\nabla_\mx (\rho_0\,V)\right]\cdot\nabla_\mx\varphi \,\md\mx  \,+\,\ds\int\rho_0\,\left( N(V)-W \right)\,\varphi\,\md\mx, \\ \, \\
W(t,\mx) \,=\,\left\{\begin{array}{l l} 0 & \text{if} \quad \rho_0(\mx) = 0,\\ \, \\
\dfrac{1}{\rho_0(\mx)}\ds\int w\,F(t,\mx,\md w) & \text{else},    \end{array}\right.
\end{array}\right.
\end{equation} 
\end{defi}

\begin{cor}\label{cor:eqlim}
Let $\Psi$ be a connectivity kernel satisfying \eqref{hyp:Psi}. Consider an initial data $(V_0,F_0)$ such that $F_0\in\mathcal{M}(\R^{d+1})$ satisfies
\begin{equation}\label{mom2}
\ds\int |w|^2\,F_0(\md\mx,\md w)\,<\,+\infty,
\end{equation}
and define for all $\mx\in\R^d$,
\begin{equation}\label{def:rho0W0}
\rho_0\,:=\,\ds\int F_0(\cdot,\md w),\quad W_0(\mx)\,:=\,\left\{\begin{array}{l l} \dfrac{1}{\rho_0(\mx)}\ds\int w\,F_0(\mx,\md w) , & \text{if } \rho_0(\mx)>0,   \\ \, \\   0, &\text{else.} \end{array}\right.
\end{equation}
Let us assume that $\rho_0$ satisfies \eqref{hyp:rho0}, $(V_0,W_0)$ satisfies \eqref{hyp:V0}, and that
\begin{equation}\label{hyp:V0=0}
V_0(\mx)\,=\,0\quad\,\text{ if }\,\mx\in\R^d\backslash\text{Supp}_\text{ess}(\rho_0).
\end{equation}
Then, for all $T>0$, there exists a unique \textcolor{black}{pair} $(V,F)$ solution of the reaction-diffusion equation \eqref{eq:lim2} on $[0,T]$ in the sense of Definition \ref{defi:sol} such that
$$\left\{\begin{array}{l}
\rho_0\,V \in L^\infty\left([0,T],H^2(\R^d)\right)\,\cap\,\mathscr{C}^0\left([0,T],H^1(\R^d)\right), \\ \, \\
F\in L^\infty\left([0,T],\mathcal{M}(\R^{d+1})\right),
\end{array}\right.$$
and such that there exists a constant $C_T>0$ such that for all $t\in[0,T]$,
$$\ds\int |w|^2\,F(t,\md\mx,\md w)\,\leq\,C_T.$$
\end{cor}

We postpone the proof of Corollary \ref{cor:eqlim} to Section \ref{sec:eqlim}. 

\subsection{Main result}
Now, we can state our main theorem about the \textcolor{black}{asymptotic limit}.
\begin{thm}\label{thm:mainresult}
Let $T>0$, and let $\Psi$ be a connectivity kernel satisfying \eqref{hyp:Psi}. Consider a set of initial data $(f_0^\ve)_{\ve>0}$ satisfying the assumptions \eqref{hyp:f0ve}-\eqref{hyp:suppf0ve}, and there exists a positive constant $M>0$ such that \begin{equation}\label{eq:ass1}
\ds\int\left( 1 + \|\mx\|^4 + |v|^4 + |w|^4 \right)\,f_0^\ve(\mx,v,w)\,\md\mx\,\md v\,\md w \,\leq\, M,
\end{equation}
\begin{equation}\label{eq:ass2}
\left\|\rho_0^\ve\right\|_{L^\infty(\R^d)} \,\leq\, M.
\end{equation}
Also consider the initial data $(\rho_0,V_0,W_0)$ satisfying the assumptions \eqref{hyp:rho0}-\eqref{hyp:V0} such that $\rho_0\in H^2(\R^d)$, and verifying:
\begin{equation}\label{eq:ass3}
\dfrac{1}{\ve^2}\|\rho_0^\ve-\rho_0\|_{L^2(\R^d)} \,\longrightarrow\,0,
\end{equation}
\begin{equation}\label{eq:ass4}
\ds\int  \rho_0^\ve(\mx) \,\left[\left|V_0^\ve(\mx)-V_0(\mx)\right|^2 + \left|W_0^\ve(\mx)-W_0(\mx)\right|^2\right]\,\md\mx \,\longrightarrow\, 0
\end{equation}
as $\ve\rightarrow0$. Consider $(V,W)$ the weak solution of the reaction-diffusion equation \eqref{eq:FHNlocal} on $[0,T]$ provided by Proposition \ref{prop:eqlim}. For any $\ve>0$, let $f^\ve$ be the weak solution of the \textcolor{black}{transport} equation \eqref{eq:f} on $[0,T]$ provided by Proposition \ref{prop:eqf}. Then, for all $\ve>0$, the macroscopic functions $(\rho_0^\ve,V^\ve,W^\ve)$ computed from $f^\ve$ satisfy
\begin{equation}\label{eq:relentrthm}
\underset{\ve\rightarrow0}{\lim}\,\underset{t\in[0,T]}{\sup}\, \ds\int \rho_0^\ve(\mx) \,\left[\left|V^\ve(t,\mx)-V(t,\mx)\right|^2 + \left|W^\ve(t,\mx)-W(t,\mx)\right|^2\right]\,\md\mx \,=\,0.
\end{equation}
Moreover, assume that there exists a measure $F_0\in\mathcal{M}(\R^{d+1})$ satisfying \eqref{mom2}-\eqref{hyp:V0=0}. Further consider $(F_0^\ve)_{\ve>0}$ the functions defined for all $\ve>0$, $(\mx,w)\in\R^{d+1}$, with
$$F_0^\ve(\mx,w)\,:=\,\ds\int_\R f_0^\ve(\mx,v,w)\,\md v.$$
Therefore, if $F_0^\ve \rightharpoonup F_0$ weakly-$\ast$ in $\mathcal{M}(\R^{d+1})$, then we have for all $\varphi\in \mathscr{C}^0_b(\R^{d+2})$,
\begin{equation}
\ds\iint \varphi(\mx,v,w)\,f^\ve(t,\mx,v,w)\,\md\mx\,\md v\,\md w \longrightarrow \int \varphi(\mx,V(t,\mx),w)\,F(t,\mx,w)\,\md\mx,
\end{equation}
strongly in $L^1_\text{loc}(0,T)$ as $\ve\rightarrow0$, where $(V,F)$ is the solution of \eqref{eq:lim} provided by Corollary \ref{cor:eqlim}.
\end{thm}

The proof is postponed to Section \ref{sec:proof}. Our approach is similar to the work done in \cite{CFF}
. To show that the \textcolor{black}{modulated energy} vanishes as $\ve$ goes to $0$, we encounter some difficulties, coming from the reaction term $\partial_v(f^\ve\,N(v))$, which makes appear some moments of $f^\ve$ of order higher than $2$ that we need to control, and from the fact that $V^\ve$ and $W^\ve$ are not \textit{a priori} differentiable in space, and not uniformly bounded. We circumvent these issues with an \textit{a priori} \textcolor{black}{dissipation} estimate, detailed in Section \ref{sec:apriori}.

\begin{rmk}
We can further precise that the convergence of the estimate \eqref{eq:relentrthm} is of order $\ve^{2/(d+6)}$ with further assumptions on the initial conditions. More precisely, we need to replace the assumption \eqref{eq:ass3} with
$$\|\rho_0^\ve-\rho_0\|_{L^2(\R^d)}\,\leq\,C\,\ve^{2+1/(d+6)},$$
to assume that the convergence of the initial conditions in \eqref{eq:ass4} is of rate $\ve^{2/(d+6)}$, and to suppose some additional regularity of $(\rho_0,V_0,W_0)$ so that $\rho_0\,V,\,\rho_0\,W \in L^\infty([0,T],H^4(\R^d)).$
\end{rmk}

\section{\textit{A priori} estimates for the \textcolor{black}{transport} equation}\label{sec:apriori}

\textcolor{black}{In this section, we prove an \textit{a priori} estimate of the moments of a solution of \eqref{eq:f}, in order to estimate a \textcolor{black}{dissipation}.} First of all, we define for all $i\in\N$ and $u\in\{v,w\}$ the moment of order $i$ in $u$ of $f^\ve$, denoted by $\mu_i^u$, and the moment of order $i$ in $\mx$ of $f^\ve$, denoted by $\mu_i^\mx$, with
$$\mu_i^u(t)\,:=\,\ds\int_{\R^{d+1}} |u|^i\,f^\ve(t,\mx,v,w)\,\md\mx\,\md v\,\md w \, , \quad \mu_i^\mx(t)\,:=\,\ds\int_{\R^{d+1}} \|\mx\|^i\,f^\ve(t,\mx,v,w)\,\md\mx\,\md v\,\md w.$$
We also define for all $\ve>0$ and for all $p\geq 1$ the \textcolor{black}{dissipation}
\begin{equation}\label{eq:defkindissip}
\mathcal{D}_p:t\,\mapsto\, \ds\dfrac{1}{2}\,\dfrac{1}{\ve^d}\iint\Psi\left(\dfrac{\|\mx-\mx'\|}{\ve}\right)\,(v^{2p-1}-v'^{2p-1})\,(v-v')\,f^\ve(t,\mz')\,f^\ve(t,\mz)\,\md\mz'\,\md\mz \,\geq\, 0,
\end{equation}
using the notation $\mz=(\mx,v,w)\in\R^{d+2}$. In the following, let $T>0$, and $\ve>0$ and suppose that there exists $f^\ve$ a well-defined solution of \eqref{eq:f} on $[0,T]$.

\begin{prop}\label{prop:moment}
Consider $f^\ve$ a weak solution of the \textcolor{black}{transport} equation \eqref{eq:f} on $[0,T]$ provided by Proposition \ref{prop:eqf}. Assume that there exists $p^\ast\in\N$ such that
\begin{equation}
\mu^v_{2p^\ast}(0)\,+\,\mu^w_{2p^\ast}(0)\,<\,+\infty.
\end{equation}
Then, for all $1\leq p\leq p^\ast$, there exists a constant $C_p$ which depends on $p$ such that for all $t\in[0,T]$, we have:
\begin{equation}\label{eq:momestimp}
\dfrac{1}{2p}\,\dfrac{\md}{\md t}\,\left(\mu^v_{2p}(t)+\mu^w_{2p}(t)\right) \,+\,\kappa_1'\, \mu^v_{2(p+1)}(t) \,+\,\dfrac{1}{\ve^2}\,\mathcal{D}_p(t) \,\leq\, C_p\,\left(\mu^v_{2p}(t)+\mu^w_{2p}(t)\right) ,
\end{equation}
where $\kappa_1'>0$ is the positive constant defined in \eqref{hyp:N}.
\end{prop}
\begin{Proof}
Let $t\in[0,T]$. In the rest of this proof, we use the notation $\bfu=(v,w)\in\R^2$ and $\mz=(\mx,\bfu)\in\R^{d+2}$. Since $f^\ve$ is a solution of the \textcolor{black}{transport} equation \eqref{eq:f}, we have:
$$\dfrac{1}{2p}\,\dfrac{\md}{\md t}\,\left(\mu_{2p}^v(t) + \mu_{2p}^w(t)\right)\,=\,I_1\,+\,I_2,$$
where
$$\left\{\begin{array}{l}
I_1\,:=\,\ds\int_{\R^{d+2}}v^{2p-1}\,\left(N(v)-w\right)\,f^\ve(t,\mz)\,\md\mz \,+\, \int_{\R^{d+2}}w^{2p-1}\,A(v,w)\,f^\ve(t,\mz)\,\md\mz, \\ \, \\
I_2\,:=\,\ds\dfrac{1}{\ve^{d+2}}\iint\Psi\left(\dfrac{\|\mx-\mx'\|}{\ve}\right)\,v^{2p-1}\,(v'-v)\,f^\ve(t,\mz')\,f^\ve(t,\mz)\,\md\mz'\,\md\mz.
\end{array}\right.$$
First of all, using Young's inequality and the properties of $N$ given in \eqref{hyp:N}, we treat the first term $I_1$ as follows:
\begin{align*}
I_1\,&\leq\,\ds\int\left( \kappa_1\,v^{2p}-\kappa_1'\,v^{2p+2}+\dfrac{2p-1}{2p}v^{2p} + \dfrac{1}{2p}w^{2p} \right)\,f^\ve(t,\mz)\,\md\mz \,   \\
&~~+\, \tau\ds\int \left( \dfrac{2p-1}{p}w^{2p} + \dfrac{1}{2p}v^{2p} - \gamma\,w^{2p} \right)\,f^\ve(t,\mz)\,\md\mz   \\
&=\, \dfrac{2p(1+\kappa_1)+\tau-1}{2p}\,\mu^v_{2p}(t)   \,+\, \dfrac{4\tau p - 2\tau + 1}{2p}\,\mu^w_{2p}(t) \,-\,\kappa_1'\,\mu^v_{2p+2}(t) .
\end{align*}
Then, to deal with the second term $I_2$, we reformulate it using the symmetry of $\Psi$. Indeed, we have:
\begin{align*}
I_2\,
&=\,-\dfrac{1}{2}\,\dfrac{1}{\ve^{d+2}}\ds\iint\Psi\left(\dfrac{\|\mx-\mx'\|}{\ve}\right)\,\left(v^{2p-1}-v'^{2p-1}\right)\,(v-v')\,f^\ve(t,\mz')\,f^\ve(t,\mz)\,\md\mz'\,\md\mz \\
&=\,-\dfrac{1}{\ve^2}\,\mathcal{D}_p(t).
\end{align*}
This enables us to conclude that there exists a constant $C_p>0$ such that for all $t\in[0,T]$,
$$
\dfrac{1}{2p}\,\dfrac{\md}{\md t}\,\left(\mu^v_{2p}(t)+\mu^w_{2p}(t)\right) \,\leq\, C_p\,\left(\mu^v_{2p}(t)+\mu^w_{2p}(t)\right) \,-\, \textcolor{black}{\kappa_1'}\mu^v_{2(p+1)}(t) \,-\,\dfrac{1}{\ve^2}\,\mathcal{D}_p(t).
$$
\end{Proof}
\begin{cor}\label{cor:moment}
Under the same assumptions than in Proposition \ref{prop:moment} with $p^\ast=2$, if we assume that for all $\ve>0$, \eqref{eq:ass1} is satisfied, then there exists a constant $C_T>0$ such that for all $k\in[0,4]$, for all $\ve>0$ and for all $t\in[0,T]$,
\begin{equation}
\left\{\begin{array}{l}
\mu_k^v(t)\,+\,\mu_k^w(t)\,+\,\mu_k^\mx(t)\,\leq\,C_T,   \\ \, \\
\ds\int_0^T\mu_{k+2}^v(t)\,\md t \,\leq\,C_T.
\end{array}\right.
\end{equation}
\end{cor}
\begin{Proof}
This result is a direct consequence of Proposition \ref{prop:moment}, integrating the inequality \eqref{eq:momestimp} between $0$ and $t$.
\end{Proof}

Now, let us estimate the \textcolor{black}{dissipation} $\mathcal{D}_1$ as defined in \eqref{eq:defkindissip}. 
\begin{cor}\label{cor:dissip0}
Let $\ve>0$. Under the same assumptions as in Proposition \ref{prop:moment} with $p^\ast=1$, there exists a constant $C_T$ such that:
\begin{equation}
\ds\int_0^T\mathcal{D}_1(t)\,\md t \,=\, \dfrac{1}{2}\,\dfrac{1}{\ve^d}\,\ds\int_0^T\iint \Psi\left(\dfrac{\|\mx-\mx'\|}{\ve}\right)\,|v-v'|^2\,f^\ve(t,\mz)\,f^\ve(t,\mz')\,\md\mz\,\md\mz'\,\md t\,\leq\,C_T\,\ve^{2}.
\end{equation}
\end{cor}
\begin{Proof}
This \textcolor{black}{dissipation} estimate comes from the inequality \eqref{eq:momestimp} with $p=1$. Indeed, integrating between $0$ and $T$, we get that 
$$\dfrac{1}{\ve^2}\,\ds\int_0^T\mathcal{D}_1(t)\,\md t\,\leq\,C\int_0^T\left( \mu_2^v(t)\,+\,\mu_2^w(t) \right)\,\md t \,+\,\mu_2^v(0)\,+\,\mu_2^w(0).$$
We conclude using the moment estimate from Corollary \ref{cor:moment}.
\end{Proof}
It turns ou that this result is not enough to conclude the proof of Theorem \ref{thm:mainresult}. 
Actually, we need to remove the weight \textcolor{black}{$\Psi_\ve(\|\cdot\|)\ast_\mx\rho_0^\ve$ in the integrand of the previous estimate. In the following, we use the shorthand notation $\Psi_\ve\ast_\mx\rho_0^\ve$}.
\color{black}
\begin{prop}\label{prop:dissip}
\color{black}
Let $\ve>0$. We make the same assumptions than in Proposition \ref{prop:moment} with $p^\ast=2$, and we further assume that there exists a function $\rho_0\in H^2(\R^d)$ such that for all $\ve>0$ small enough,
\begin{equation}\label{eq:hyprho0}
\|\rho_0^\ve \,-\,\rho_0\|_{L^2(\R^d)}\,\leq\, C\,\ve^2,
\end{equation}
for some positive constant $C>0$. Then, there exists a constant $C_T$ such that for all $t\in[0,T]$, we have:
\begin{equation}
\ds\int_0^T\int_{\R^{d+2}}f^\ve(t,\mx,v,w)\,\left| v - V^\ve(t,\mx) \right|^2\,\md\mx\,\md v\,\md w\,\leq\,C_T\,\ve^{4/(d+6)}.
\end{equation}
\end{prop}
\color{black}
\begin{Proof}
Let $\ve>0$ and $T>0$. We define the integral
$$I_\ve\,:=\,\ds\int_0^T\int_{\R^{d+2}}f^\ve(t,\mx,v,w)\,\left| v - V^\ve(t,\mx) \right|^2\,\md\mx\,\md v\,\md w\,\md t.$$
Using the definition of $V^\ve$, let us notice that {\color{black}for all $t\in[0,T]$ and all $\mx\in\R^d$,
$$\ds\int_{\R^2}f^\ve(t,\mx,v,w)\,\left| v - V^\ve(t,\mx) \right|^2\,\md v\,\md w \,=\,  \ds\int_{\R^2}f^\ve(t,\mx,v,w)\,|v|^2\,\md v\,\md w \,-\, \rho_0^\ve(\mx)\,\left| V^\ve(t,\mx) \right|^2.$$}
In the rest of this proof, we use the notation $\mz=(\mx,v,w)\in\R^{d+2}$. First of all, \textcolor{black}{we restrict our analysis to the nontrivial case $\rho_0^\ve\not\equiv 0$. Since $\Psi>0$ almost everywhere and $\rho_0^\ve\geq 0$, we get that for all $\mx\in\R^d$, $\Psi_\ve\ast_\mx\rho_0^\ve(\mx) \,>\,0.$}
Our strategy to estimate $I_\ve$ consists in dividing the set of integration into subsets on which the integrand is easier to control. Let $\eta>0$ be a constant depending on $\ve$ to be determined later. We define:
$$\left\{\begin{array}{l}
\mathcal{A}^\eta_\ve \,:=\, \left\{ \mx\in\R^d \, , \, \Psi_\ve\ast_\mx\rho_0^\ve(\mx) \,\geq\, \eta \right\}, \\ \, \\
\mathcal{B}^\eta_\ve \,:=\, \left\{ \mx\in\R^d \, , \, 0 \,<\,\Psi_\ve\ast_\mx\rho_0^\ve(\mx) \,<\, \eta \right\},
\end{array}\right.$$
and hence $\R^d = \mathcal{A}^\eta_\ve \cup \mathcal{B}^\eta_\ve$. Then, we have that
\begin{align*}
\ds\int_0^T\underset{\mathcal{A}^\eta_\ve\times\R^2}{\iint}f^\ve\,&|v-V^\ve|^2\,\md\mz\,\md t \,\leq\,\dfrac{1}{\eta}\ds\int_0^T\underset{\mathcal{A}^\eta_\ve\times\R^2}{\iint}f^\ve\,|v-V^\ve|^2\,\Psi_\ve\ast_\mx\rho_0^\ve(\mx)\,\md\mz\,\md t  \\
&=\,\dfrac{1}{\eta}\ds\int_0^T\underset{\mathcal{A}^\eta_\ve\times\R^2}{\iint}f^\ve\,|v|^2\,\Psi_\ve\ast_\mx\rho_0^\ve(\mx)\,\md\mz\,\md t   \,-\, \dfrac{1}{\eta}\ds\int_0^T\underset{\mathcal{A}^\eta_\ve}{\int}\rho_0^\ve\,|V^\ve|^2\,\Psi_\ve\ast_\mx\rho_0^\ve(\mx)\,\md\mx\,\md t    \\
&=\,\dfrac{1}{\eta}\ds\int_0^T\underset{\mathcal{A}^\eta_\ve\times\R^2}{\iint}f^\ve\,|v|^2\,\Psi_\ve\ast_\mx\rho_0^\ve(\mx)\,\md\mz\,\md t   \,-\, \dfrac{1}{\eta}\ds\int_0^T\underset{\mathcal{A}^\eta_\ve}{\int}\rho_0^\ve\,V^\ve\,\Psi_\ve\ast_\mx[\rho_0^\ve\,V^\ve](\mx)\,\md\mx\,\md t    \\
&~~~+\,\dfrac{1}{\eta}\ds\int_0^T\underset{\mathcal{A}^\eta_\ve\times\R^2}{\iint}\rho_0^\ve\,V^\ve\,\Psi_\ve\ast_\mx[\rho_0^\ve\,V^\ve](\mx)\,\md\mz\,\md t   \,-\, \dfrac{1}{\eta}\ds\int_0^T\underset{\mathcal{A}^\eta_\ve}{\int}\rho_0^\ve\,|V^\ve|^2\,\Psi_\ve\ast_\mx\rho_0^\ve(\mx)\,\md\mx\,\md t    \\
&=\,\dfrac{1}{\eta}\dfrac{1}{\ve^d}\dfrac{1}{2}\int_0^T\iint\Psi\left(\dfrac{\|\mx-\mx'\|}{\ve}\right)\,|v-v'|^2\,f^\ve(t,\mz)\,f^\ve(t,\mz')\,\md\mz'\,\md\mz\,\md t    \\
&~~~-\,\dfrac{1}{\eta}\dfrac{1}{\ve^d}\dfrac{1}{2}\int_0^T\iint\Psi\left(\dfrac{\|\mx-\mx'\|}{\ve}\right)\,|V^\ve(t,\mx)-V^\ve(t,\mx')|^2\,\rho_0^\ve(\mx)\,\rho_0^\ve(\mx')\,\md\mx'\,\md\mx\,\md t   \\
&\leq\,\dfrac{1}{\eta}\,\ds\int_0^T \mathcal{D}_1(t) \,\md t,
\end{align*}
where $\mathcal{D}_1$ is defined in \eqref{eq:defkindissip}. Consequently, using Corollary \ref{cor:dissip0}, we conclude that there exists a positive constant $C_T>0$ such that
\begin{equation}\label{eq:Aeps}
\ds\int_0^T\underset{\mathcal{A}^\eta_\ve\times\R^2}{\iint}f^\ve\,|v-V^\ve|^2\,\md\mz\,\md t \,\leq\, C_T\,\dfrac{\ve^2}{\eta}.
\end{equation}
Then, it remains to estimate
$$\,\ds\int_0^T\underset{\mathcal{B}^\eta_\ve\times\R^2}{\iint}f^\ve\,|v-V^\ve|^2\,\md\mz\,\md t \,\leq\, \ds\int_0^T\underset{\mathcal{B}^\eta_\ve\times\R^2}{\iint}f^\ve\,|v|^2\,\md\mz\,\md t \,=\, I_1+I_2+I_3,$$
where 
$$\left\{\begin{array}{l}
I_1 \,:=\, \ds\int_0^T\int_{\mathcal{B}^\eta_\ve}\int_{\{|v|> R\}}f^\ve\,|v|^2\,\md\mz\,\md t, \\ \, \\
I_2 \,:=\, \ds\int_0^T\int_{\mathcal{B}^\eta_\ve\cap B^c(0,R)}\int_{\{|v|\leq R\}}f^\ve\,|v|^2\,\md\mz\,\md t, \\ \, \\
I_3 \,:=\, \ds\int_0^T\int_{\mathcal{B}^\eta_\ve\cap B(0,R)}\int_{\{|v|\leq R\}}f^\ve\,|v|^2\,\md\mz\,\md t,
\end{array}\right.$$
where $R>0$ is a constant depending on $\eta$ and $\ve$ to be determined later. For $k>2$, we get that
$$I_1\,\leq\,\dfrac{1}{R^{k-2}}\,\ds\int_0^T\int_{\mathcal{B}^\eta_\ve}\int_{\{|v|> R\}} f^\ve\,|v|^k\,\md \mz\,\md t \,\leq\,\dfrac{1}{R^{k-2}}\,\ds\int_0^T\mu_k^v(t)\,\md t .$$
Then, for $q>2$, we also have that
$$I_2\,\leq\,\ds\int_0^T\int f^\ve\,R^2\,\dfrac{|\mx|^q}{R^q}\,\md\mz\,\md t\,\leq\,\dfrac{1}{R^{q-2}}\,\ds\int_0^T\mu_q^\mx(t)\,\md t .$$
As for the last term, we compute:
\begin{align*}
I_3\,&\leq\,R^2\ds\int_0^T\int_{\mathcal{B}^\eta_\ve\cap B(0,R)}\rho_0^\ve(\mx)\,\md\mx\,\md t    \\
&\leq\,R^2\ds\int_0^T\int_{\mathcal{B}^\eta_\ve\cap B(0,R)} \Psi_\ve\ast\rho_0^\ve(\mx) \,\md\mx\,\md t \,+\, R^2\int_0^T\int_{\mathcal{B}^\eta_\ve\cap B(0,R)} \left|\rho_0^\ve(\mx) \,-\, \Psi_\ve\ast\rho_0^\ve(\mx)\right| \,\md\mx\,\md t    \\
&\leq\,C\,T\,R^{d+2}\,\eta \,+\, R^2\,T\,\left( \ds\int_{\mathcal{B}^\eta_\ve\cap B(0,R)}1\,\md\mx \right)^{1/2}\,\|\rho_0^\ve - \Psi_\ve\ast\rho_0^\ve\|_{L^2(\R^d)}  \\
&\leq\,C\,T\,R^{d+2}\,\eta \,+\, C^{1/2}\,T\,R^{2+d/2}\|\rho_0^\ve - \Psi_\ve\ast\rho_0^\ve\|_{L^2(\R^d)}.
\end{align*}
Then, using Young's inequality, we notice that
\begin{align*}
\|\rho_0^\ve - \Psi_\ve\ast\rho_0^\ve\|_{L^2(\R^d)} \,&\leq\, \|\rho_0^\ve - \rho_0\|_{L^2(\R^d)} \,+\, \|\rho_0 - \Psi_\ve\ast\rho_0\|_{L^2(\R^d)} \,+\, \|\Psi_\ve\ast(\rho_0^\ve-\rho_0)\|_{L^2(\R^d)}   \\
&\leq\,(1+\|\Psi_\ve\|_{L^1(\R^d)})\,\|\rho_0^\ve-\rho_0\|_{L^2(\R^d)} \,+\, \|\rho_0 - \Psi_\ve\ast\rho_0\|_{L^2(\R^d)}   \\
&\leq\,2\,\|\rho_0^\ve-\rho_0\|_{L^2(\R^d)} \,+\, \|\rho_0 - \Psi_\ve\ast\rho_0\|_{L^2(\R^d)}.
\end{align*}
On the one hand, we have assumed that the initial data satisfies the estimate \eqref{eq:hyprho0}. On the other hand, we can estimate $\|\rho_0 - \Psi_\ve\ast\rho_0\|_{L^2(\R^d)}$ with similar arguments as in \cite{BATES}. Indeed, using the change of variable $\my=(\mx-\mx')/\ve$ and a Taylor expansion, we get that
\begin{align*}
\|\rho_0 - \Psi_\ve\ast\rho_0\|_{L^2(\R^d)}^2\,&\leq\,\ds\int\left| \dfrac{1}{\ve^d}\int\Psi\left( \dfrac{\mx-\mx'}{\ve} \right)\,(\rho_0(\mx')-\rho_0(\mx)) \,\md\mx'\right|^2\,\md\mx   \\
&\leq\,\ds\int\left| \int\Psi\left( \my \right)\,(\rho_0(\mx-\ve\my)-\rho_0(\mx)) \,\md\my\right|^2\,\md\mx   \\
&\leq\,\ve^4\ds\int\left| \int\Psi\left( \my \right)\,\int_0^1 (1-s)\,\my^T\cdot\nabla_\mx^2\rho_0(\mx-\ve s \my)\cdot\my\,\md s \,\md\my\right|^2\,\md\mx   .
\end{align*}
Furthermore, using Cauchy-Schwarz inequality for the integral in $s$ and then for the integral in $\my$, we get
\begin{align*}
\|\rho_0 - &\Psi_\ve\ast\rho_0\|_{L^2(\R^d)}^2\,   \\
&\leq\,\ve^4\ds\int\left| \int\Psi\left( \my \right)\left( \int_0^1|1-s|\|\my\|^2\,\md s \right)^{1/2}\left( \int_0^1|1-s|\|\my\|^2\,\|\nabla_\mx^2\rho_0(\mx-\ve s \my) \|^2\md s \right)^{1/2}\md\my\right|^2\,\md\mx   \\
&\leq\,\ve^4\ds\int\left(\int\Psi(\|\my\|)\dfrac{\|\my\|^2}{2}\,\md\my\right)\,\left( \int\Psi(\|\my\|)\int_0^1|1-s|\|\my\|^2\,\|\nabla_\mx^2\rho_0(\mx-\ve s \my) \|^2\,\md s\,\md\my \right)\,\md\mx
\end{align*}
Consequently, 
\begin{align*}
\|\rho_0 - \Psi_\ve\ast\rho_0\|_{L^2(\R^d)}^2 \,&\leq\, \sigma\,\ve^4\ds\int\Psi(\|\my\|)\int_0^1|1-s|\|\my\|^2\int \|\nabla_\mx^2\rho_0(\mx-\ve s \my) \|^2 \,\md\mx\,\md s\,\md \my   \\
&\leq\,\sigma^2\,\|\rho_0\|_{H^2(\R^d)}^2\,\ve^4.
\end{align*}
Finally, we get that there exists a positive constant $C>0$ such that 
\begin{equation}\label{eq:I3}
I_3\,\leq\,C\,T\,\left( R^{d+2}\,\eta \,+\, R^{(d+4)/2}\,\ve^2 \right).
\end{equation}
Finally, using the moment estimates, and the estimate \eqref{eq:Aeps}, we get that there exists a positive constant $C_T$ such that
\begin{equation}
I_\ve \,\leq\, C_T\,\left( \dfrac{\ve^2}{\eta} \,+\, R^{d+2}\,\eta \,+\, R^{(d+4)/2}\,\ve^2 \,+\, \dfrac{1}{R^{k-2}} \,+\,\dfrac{1}{R^{q-2}} \right).
\end{equation}
It remains to optimize the values of $\eta$ and $R$. For the sake of simplicity, we choose $k=q=4$. We consider $R = \eta^{-1/(d+4)}$, so that 
$$R^{d+2}\eta \,=\,\dfrac{1}{R^2}.$$
Then, we take $\eta = \ve^{2(d+4)/(d+6)}$, so that
$$R^{d+2}\eta\,=\,\dfrac{\ve^2}{\eta}\textcolor{black}{\,=\,\ve^{4/(d+6)}}.$$
This leads to
\begin{equation}
I_\ve \,\leq\, C_T\,\left( \ve^{4/(d+6)} \,+\, \ve^{(d+8)/(d+6)} \right) \,\leq\, \widetilde{C}_T \,  \ve^{4/(d+6)},
\end{equation}
for a positive constant $\widetilde{C}_T>0$ and $\ve>0$ small enough.
\end{Proof}

\section{Proof of Theorem \ref{thm:mainresult}}\label{sec:proof}

{\color{black}
Our proof of Theorem \ref{thm:mainresult} relies on a modulated energy argument, as developed in \cite{BRE}. This leads to estimate the distance between the macroscopic functions derived from the solution of the \textcolor{black}{transport} equation \eqref{eq:f}, and the solution of the limit system \eqref{eq:lim}. First, we introduce the notion of modulated energy we use in this article. Then, we prove that it converges to $0$ as $\ve$ goes to $0$. Finally, we explain how this argument enables us to prove Theorem \ref{thm:mainresult}.
}

In the rest of this article, for any given $\ve>0$ and for $\rho:\R^d\rightarrow\R$ and $V:(0,\infty)\times\R^d\rightarrow\R$ regular enough, we define the following local and nonlocal differential operators:
\begin{equation}\label{eq:op}
\left\{\begin{array}{l}
\mathcal{L}_\rho(V) \,:=\,\sigma \big[ \Delta_\mx(\rho\,V) \,-\, \Delta_\mx\rho \, V \big] \,=\, \sigma \big[ \rho\,\Delta_\mx V \,+\, 2\,\nabla_\mx\rho\cdot\nabla_\mx V \big], \\ \, \\
\mathscr{L}_\rho(V) \,:=\, \Psi_\ve\ast_\mx\left[\rho\,V\right] \,-\,\left[\Psi_\ve\ast_\mx\rho\right]\,V \,=\, \dfrac{1}{\ve^{d+2}}\ds\iint\Psi\left(\dfrac{\|\mx-\mx'\|}{\ve}\right)\,\left( V(t,\mx')-V(t,\mx) \right)\,\rho(\mx')\,\md\mx',
\end{array}\right.
\end{equation}
respectively defined on $H^2(\R^d)$ and $L^\infty(\R^d)$.
\subsection{Definition of \textcolor{black}{modulated energy}}\label{sec:def}

{\color{black}
Consider $\mZ^\ve=(\rho_0^\ve,\rho_0^\ve\,V^\ve,\rho_0^\ve\,W^\ve)$ the triple of macroscopic quantities computes from $f^\ve$ the solution to the transport equation \eqref{eq:f}, and $\mZ=(\rho_0,\rho_0\,V,\rho_0\,W)$ the solution of the reaction-diffusion equation \eqref{eq:rhoFHNlocal}. Then, we define the modulated energy of our system as follows for all $t>0$:
\begin{equation}
\mathcal{H}_\ve(t)\,:=\,\ds\int_{\R^d}\rho_0^\ve\,\dfrac{|V-V^\ve|^2\,+\,|W-W^\ve|^2}{2}\,\md\mx.
\end{equation}
}

\subsection{\textcolor{black}{Modulated energy} estimate}\label{subsec:estim}

This subsection is devoted to the proof of the \textcolor{black}{modulated energy} estimate \eqref{eq:relentrthm} under the same assumptions as in Theorem \ref{thm:mainresult}. 
For all $\ve>0$, let $f^\ve$ be the solution of the \textcolor{black}{transport} equation \eqref{eq:f}. According to Corollary \ref{cor:moment}, we know that for all $t\in[0,T]$, the moment of order $4$ of $f^\ve(t)$ is uniformly bounded with respect to $\ve>0$. Therefore, using Hölder's inequality, we obtain that for all $\mx\in\R^d$ such that $\rho_0^\ve(\mx)>0$ and for all $t\in[0,T]$, 
\begin{eqnarray}
\rho_0^\ve(\mx)\,|V^\ve(t,\mx)|^4 &=& \dfrac{1}{|\rho_0^\ve(\mx)|^3}\left( \ds\int v\,f^\ve(t,\mx,v,w)\,\md v\,\md w  \right)^4   \nonumber  \\
&\leq & \ds\int |v|^4\,f^\ve(t,\mx,v,w)\,\md v\md w.\label{eq:rhoV4}
\end{eqnarray}
This last inequality \eqref{eq:rhoV4} remains true where $\rho_0^\ve(\mx)\,=\,0$ and with $W^\ve$ instead of $V^\ve$. Consequently, since $\|\rho_0^\ve\|_{L^1(\R^d)}=1$, we get that for any $0\leq p\leq 4$ and for all $t\in[0,T]$, 
$$\rho_0^\ve\left( |V^\ve(t)|^p + |W^\ve(t)|^p \right)\in L^1(\R^d).$$
{\color{black}
Let $(V,W)$ be the weak solution of the reaction-diffusion system \eqref{eq:FHNlocal} provided by Proposition \ref{prop:eqlim}. In the following, we consider the triples $\mZ^\ve$ and $\mZ$ as defined in Subsection \ref{sec:def}.} Since $V$ and $W$ are in $W^{1,\infty}([0,T],L^2(\R^d))$ by definition, for all $t\in[0,T]$, we can compute :
\begin{align*}
\mathcal{H}_\ve(t)\,&=\,\mathcal{H}_\ve(0)  \\
&~~\,+\,\ds\int_0^t\left[\int\left( V^\ve-V \right)\,\left( \partial_t(\rho_0^\ve\,V^\ve) \,-\, \rho_0^\ve\partial_t V \right)\,\md\mx \,+\, \ds\int\left( W^\ve-W \right)\,\left( \partial_t(\rho_0^\ve\,W^\ve) \,-\, \rho_0^\ve\partial_t W \right)\,\md\mx \right](s)\,\md s  \\
&=\,\mathcal{H}_\ve(0) \,+\, \ds\int_0^t\left[\mathcal{T}_{1}(s) \,+\, \mathcal{T}_{2}(s) \,+\,\mathcal{T}_{3}(s)\right]\,\md s ,
\end{align*}
where for all $s\in[0,T]$, we define
$$\left\{\begin{array}{l}
\mathcal{T}_1(s) \,:=\, \ds\int\rho_0^\ve\,(W^\ve-W)\,\left(A(V^\ve,W^\ve)\,-\, A(V,W)\right)\,\md\mx \,-\,\ds\int\rho_0^\ve\,(V^\ve-V)\,\left(W^\ve-W\right)\,\md\mx, \\ \, \\
\mathcal{T}_2(s) \,:=\, \ds\int_{\R^d}(V^\ve-V)\,\ds\int_{\R^2} \left(N(v)-N(V)\right)\,f^\ve(s,\mx,\bfu)\,\md\bfu\,\md\mx, \\ \, \\
\mathcal{T}_3(s) \,:=\, \ds\int\rho_0^\ve\,(V^\ve-V)\,\left(\mathscr{L}_{\rho_0^\ve}(V^\ve) \,-\,\mathcal{L}_{\rho_0}(V)\right)\,\md\mx,
\end{array}\right.$$
which respectively stand for the difference between the linear reaction terms, the nonlinear reaction terms, and the diffusion terms. 
\paragraph{Estimate of the linear reaction terms. } First of all, we can directly treat the first term $\mathcal{T}_1$ with Young's inequality, which yields that,
\begin{equation}\label{eq:T1}
\ds\int_0^T\mathcal{T}_1(t)\,\md t \,\leq\,(1+\tau)\,\ds\int_0^T\textcolor{black}{\mathcal{H}_\ve(t)},\md t.
\end{equation}
\paragraph{Estimate of the nonlinear reaction terms. } Then, we deal with the second term $\mathcal{T}_2$ using the assumptions \eqref{hyp:N} satisfied by $N$, as \textcolor{black}{in \cite{CFF}. For all $t\in[0,T]$, we have:} 
\begin{align*}
\mathcal{T}_2(t)\,&=\,\ds\int_{\R^d}(V^\ve-V)\,\ds\int_{\R^2} \left(N(v)-N(V^\ve)\right)\,f^\ve(t,\mx,\bfu)\,\md\bfu\,\md\mx \,+\, \ds\int_{\R^d}(V^\ve-V)\,\left(N(V^\ve)-N(V)\right)\,\rho_0^\ve(\mx)\,\md\mx   \\
&\leq\,\kappa_3\ds\int_{\R^{d+2}}|V^\ve-V|\,|V^\ve-v|\,\left[1+v^2+(V^\ve)^2\right]\,f^\ve(t,\mx,\bfu)\,\md\mx\,\md\bfu \,+\,2\,\kappa_2\,\textcolor{black}{\mathcal{H}_\ve(t)},
\end{align*}
where the constants $\kappa_2$ and $\kappa_3$ are given in \eqref{hyp:N}. Then, in order to estimate $\mathcal{T}_2$ using the \textcolor{black}{dissipation} estimate from \textcolor{black}{Proposition \ref{prop:dissip}}, Cauchy-Schwarz inequality yields that
$$\mathcal{T}_2(t)\,\leq\,\alpha(t)\,\left( \ds\int \left|V^\ve(t)-v\right|^2\,f^\ve(t,\mx,\bfu)\,\md\mx\,\md\bfu \right)^{1/2} \,+\,2\,\kappa_2\,\textcolor{black}{\mathcal{H}_\ve(t)},$$
where
$$\alpha(t)\,:=\,\kappa_3\,\left( \ds\int\left[ 1 + (V^\ve(t))^2 + v^2 \right]^2\,\left[ V^\ve(t) - V(t) \right]^2\,f^\ve(t,\mx,\bfu)\,\md\mx\,\md\bfu \right)^{1/2}.$$
We recall that $V\in L^\infty([0,T],H^2(\R^d))$, and $H^2(\R^d)\,\subset\, L^\infty(\R^d)$ since $d\leq3$. Hence, using the moment estimate from Corollary \ref{cor:moment}, and the fact that for all $t\in[0,T]$ and $\mx\in\R^d$,
$$\rho_0^\ve(\mx)\,|V^\ve(t,\mx)|^6\,\leq\,\ds\int |v|^6\,f^\ve(t,\mx,\bfu)\,\md\mx\,\md\bfu,$$
we can conclude that there exists a positive constant $C_T>0$ such that 
$$\ds\int_0^T\alpha(t)^2\,\md t\,\leq\, C_T.$$
Consequently, according to the estimate from \textcolor{black}{Proposition \ref{prop:dissip}}, 
\begin{equation}\label{eq:T2}
\ds\int_0^T \mathcal{T}_2(t)\,\md t \,\leq\, C_T\,\ve^{2/(d+6)}\,+\,2\,\kappa_2\ds\int_0^T\textcolor{black}{\mathcal{H}_\ve(t)}\,\md t.
\end{equation}
\paragraph{Estimate of the diffusion terms. } Finally, it remains to estimate the third term $\mathcal{T}_3$, involving the difference between the nonlocal diffusion term $\mathscr{L}_{\rho_0^\ve}(V^\ve)$ and the local diffusion term $\mathcal{L}_{\rho}(V)$. On the one hand, since $V^\ve$ is not regular enough in space, we cannot apply the operator $\mathcal{L}_{\rho_0}$ to it. On the other hand, we can apply the nonlocal operator $\mathscr{L}_{\rho_0^\ve}$ to both $V^\ve$ and $V$. This leads to rewrite $\mathcal{T}_3$ as follows:
$$\mathcal{T}_3\,=\,\mathcal{T}_{3,1}\,+\, \mathcal{T}_{3,2},$$
where for all $t\in[0,T]$
$$\left\{\begin{array}{l}
\mathcal{T}_{3,1}(t)\,:=\,\ds\int\rho_0^\ve\,\left( V^\ve-V \right)\,\left( \mathscr{L}_{\rho_0^\ve}(V^\ve) \,-\, \mathscr{L}_{\rho_0}(V) \right)\,\md\mx,  \\ \, \\
\mathcal{T}_{3,2}(t)\,:=\,\ds\int\rho_0^\ve\,\left( V^\ve-V \right)\,\left( \mathscr{L}_{\rho_0}(V) \,-\, \mathcal{L}_{\rho_0}(V) \right)\,\md\mx.
\end{array}\right.$$
To estimate the first term $\mathcal{T}_{3,1}$, using the shorthand notations $V:=V(t,\mx)$, $V':=V(t,\mx')$, and the same for $V^\ve$ and $V^\ve\,'$, we compute:
\begin{align*}
\mathcal{T}_{3,1}(t)\,&=\,\dfrac{1}{\ve^{d+2}}\ds\iint\Psi\left( \dfrac{\|\mx-\mx'\|}{\ve} \right)\,\rho_0^\ve(\mx)\,\left(V^\ve-V\right)\left[\rho_0^\ve(\mx')\,\left(V^\ve\ds'-V^\ve\right) \,-\, \rho_0(\mx')\,\left( V'-V\right) \right]\,\md\mx\,\md\mx'   \\
&=\,\dfrac{1}{\ve^{d+2}}\ds\iint\Psi\left(\dfrac{\|\mx-\mx'\|}{\ve}\right)\,\left(V^\ve-V\right)\,\rho_0^\ve(\mx)\rho_0^\ve(\mx')\left[  (V^\ve\,'-V') \,-\, (V^\ve-V)  \right]\,\md\mx\,\md\mx'   \\
&~~~+\,\dfrac{1}{\ve^{d+2}}\ds\iint\Psi\left(\dfrac{\|\mx-\mx'\|}{\ve}\right)\,\left(V^\ve-V\right)\,\rho_0^\ve(\mx)\,(\rho_0^\ve(\mx')-\rho_0(\mx'))\,(V'-V)\,\md\mx\,\md\mx'   \\
&\leq\,\dfrac{1}{\ve^{d+2}}\ds\iint\Psi\left(\dfrac{\|\mx-\mx'\|}{\ve}\right)\,\rho_0^\ve(\mx)\,\left| V^\ve-V \right|\,\left| \rho_0^\ve-\rho_0 \right|(\mx')\,\left|V'-V\right|\,\md\mx\,\md\mx'   \\
&\leq\,\,2\,\|V\|_{L^\infty}\,\dfrac{1}{\ve^{d+2}}\iint\Psi\left(\dfrac{\|\mx-\mx'\|}{\ve}\right)\,\rho_0^\ve(\mx)\,\left| V^\ve-V \right|\,\left| \rho_0^\ve-\rho_0 \right|(\mx')\,\md\mx\,\md\mx',
\end{align*}
and then, using Young's inequality, we have
\begin{align*}
\mathcal{T}_{3,1}(t) \,&\leq\,\,\|V\|_{L^\infty}\, \ds\int\rho_0^\ve\,  \left| V-V^\ve \right|^2\,\md\mx \,+\,\|V\|_{L^\infty}\,\int\left[\dfrac{1}{\ve^{d+2}}\int\Psi\left(\dfrac{\|\mx-\mx'\|}{\ve}\right)\,\left| \rho_0^\ve-\rho_0 \right|(\mx')\,\md\mx'\right]^2\,\rho_0^\ve(\mx)\,\md\mx   \\
&\leq\,2\,\|V\|_{L^\infty}\, \textcolor{black}{\mathcal{H}_\ve(t)} \,+\,\dfrac{1}{\ve^4}\,\|V\|_{L^\infty}\,\|\rho_0^\ve\|_{L^\infty}\,\|\Psi_\ve\|_{L^1}\,\|\rho_0-\rho_0^\ve\|_{L^2}^2.
\end{align*}
This leads to the estimate
\begin{equation}\label{eq:T31}
\mathcal{T}_{3,1}(t)\,\leq\, C_T\left( \dfrac{1}{\ve^4}\,\|\rho_0^\ve-\rho_0\|_{L^2(\R^d)}^2 \,+\,\textcolor{black}{\mathcal{H}_\ve(t)} \right),
\end{equation}
where $C_T>0$ is a positive constant independent of $\ve$. It remains to control the final term $\mathcal{T}_{3,2}$. We start by separating the diffusions on $\rho_0\,V$ and on $\rho_0$ alone, as follows:
$$\mathcal{T}_{3,2}\,=\,\mathcal{T}_{3,2,1}\,+\,\mathcal{T}_{3,2,2},$$
where for all $t\in[0,T]$
$$\left\{\begin{array}{l}
\mathcal{T}_{3,2,1}(t)\,:=\, \ds\int\rho_0^\ve\,(V^\ve-V)\,\left[ \dfrac{1}{\ve^{2}}\left(\Psi_\ve\ast_\mx[\rho_0\,V](t,\mx) \,-\,\rho_0\,V(t,\mx)\right)\,-\,\sigma\,\Delta_\mx (\rho_0\,V)(t,\mx) \right]\,\md\mx, \\ \, \\
\mathcal{T}_{3,2,2}(t)\,:=\, -\ds\int\rho_0^\ve\,(V^\ve-V)\,V(t,\mx)\left[ \dfrac{1}{\ve^{2}}\left(\Psi_\ve\ast_\mx\rho_0(\mx) \,-\,\rho_0(\mx)\right)\,-\,\sigma\,\Delta_\mx \rho_0(\mx) \right]\,\md\mx.
\end{array}\right.$$
Our strategy to estimate both $\mathcal{T}_{3,2,1}$ and $\mathcal{T}_{3,2,2}$ follows the idea from \cite{BATES} with a Taylor expansion. Using Young's inequality, we get that
$$\left\{\begin{array}{l}
\mathcal{T}_{3,2,1}(t)\,\leq\, \textcolor{black}{\mathcal{H}_\ve(t)} \,+\,\dfrac{1}{2}\,\|\rho_0^\ve\|_{L^\infty}\ds\int\left| \dfrac{1}{\ve^2}\left( \Psi_\ve\ast_\mx[\rho_0\,V](t,\mx)\,-\,\rho_0\,V(t,\mx) \right) \,-\,\sigma\,\Delta_\mx(\rho_0\,V)(t,\mx) \right|^2\,\md\mx,  \\ \, \\
\mathcal{T}_{3,2,2}(t)\,\leq\, \textcolor{black}{\mathcal{H}_\ve(t)} \,+\,\dfrac{1}{2}\|\rho_0^\ve\|_{L^\infty}\,\|V\|_{L^\infty}^2\ds\int\left| \dfrac{1}{\ve^2}\left( \Psi_\ve\ast_\mx\rho_0(\mx)\,-\,\rho_0(\mx) \right) \,-\,\sigma\,\Delta_\mx\rho_0(\mx) \right|^2\,\md\mx.
\end{array}\right.$$
Then, for all $\mx\in\R^d$, we apply in the convolution products 
the change of variable $\my=(\mx-\mx')/\ve$, so that using a Taylor expansion, we get:
\begin{align*}
&\left| \dfrac{1}{\ve^2}\left( \Psi_\ve\ast_\mx[\rho_0\,V](t,\mx)\,-\,\rho_0\,V(t,\mx) \right) \,-\,\sigma\,\Delta_\mx(\rho_0\,V)(t,\mx) \right|^2   \\
&=\,\left| \dfrac{1}{\ve^2}\ds\int\Psi(\|\my\|)\,(\rho_0\,V)(t,\mx-\ve\my)\,\md\my\,-\dfrac{1}{\ve^2}\,\rho_0\,V(t,\mx)\,-\,\sigma\,\Delta_\mx(\rho_0\,V)(t,\mx)  \right|^2   \\
&=\,\left|\ds\int\Psi(\|\my\|)\ds\int_0^1 (1-s)\,\my^T\cdot\left( \nabla_\mx^2(\rho_0\,V)(t,\mx-\ve\,s\,\my) - \nabla_\mx^2(\rho_0\,V)(t,\mx) \right) \cdot\my\,\md s\,\md\my \right|^2.
\end{align*}
Besides, we consecutively use the Cauchy-Schwarz inequality in the integrals in $s$ and then in $\my$, which gives:
\begin{align*}
&\left| \dfrac{1}{\ve^2}\left( \Psi_\ve\ast_\mx[\rho_0\,V](t,\mx)\,-\,\rho_0\,V(t,\mx) \right) \,-\,\sigma\,\Delta_\mx(\rho_0\,V)(t,\mx) \right|^2  \\
&\leq\, \left| \ds\int\Psi(\|\my\|)\left(\ds\int_0^1 |1-s|\,\|\my\|\,\md s  \right)^{1/2}\left( \ds\int_0^1 |1-s|\,\|\my\|^2\,\left\|  \nabla_\mx^2(\rho_0\,V)(t,\mx-\ve\,s\,\my) - \nabla_\mx^2(\rho_0\,V)(t,\mx)  \right\|^2  \,\md s \right)^{1/2}   \md\my\right|^2    \\
&\leq\,\sigma\ds\int \Psi(\|\my\|)\,\left( \ds\int_0^1 |1-s|\,\|\my\|^2\,\left\|  \nabla_\mx^2(\rho_0\,V)(t,\mx-\ve\,s\,\my) - \nabla_\mx^2(\rho_0\,V)(t,\mx)  \right\|^2  \,\md s \right)   \,\md\my .
\end{align*}
Consequently, after integrating this last inequality with respect to $\mx$, using the fact that $\rho_0\,V$ is in $L^\infty([0,T],H^2(\R^d))$ and the hypotheses \eqref{hyp:Psi} satisfied by $\Psi$, we get:
\begin{align*}
\ds\int\bigg| \dfrac{1}{\ve^2}&\left( \Psi_\ve\ast_\mx[\rho_0\,V](t,\mx)\,-\,\rho_0\,V(t,\mx) \right) \,-\,\sigma\,\Delta_\mx(\rho_0\,V)(t,\mx) \bigg|^2 \,\md\mx  \\
&\leq\,\sigma\ds\iint \Psi(\|\my\|)\,\left( \ds\int_0^1 |1-s|\,\|\my\|^2\,\left\|  \nabla_\mx^2(\rho_0\,V)(t,\mx-\ve\,s\,\my) - \nabla_\mx^2(\rho_0\,V)(t,\mx)  \right\|^2  \,\md s \right)   \,\md\my\,\md\mx \\
&\leq\,2\,\sigma^2\,\ds\|\rho_0\,V\|_{L^\infty([0,T],H^2(\R^d))}^2.
\end{align*}
Then, using these two last inequalities and the Lebesgue's Dominated Convergence Theorem, and the fact that $\|\rho_0^\ve\|_{L^\infty}$ is uniformly bounded, we get that as $\ve$ goes to $0$, for all $t\in[0,T]$
\begin{equation}\label{eq:T321}
\mathcal{T}_{3,2,1}(t)\,\leq\,\textcolor{black}{\mathcal{H}_\ve(t)}\,+\,o_{\ve\rightarrow 0}(1),
\end{equation}
where $o_{\ve\rightarrow 0}(1)$ denotes a function which converges towards $0$ as $\ve$ goes to $0$, uniformly in $\textcolor{black}{t\in[0,T]}$. Using similar arguments, and the fact that $\rho_0\in H^2(\R^d)$ and $V\in L^\infty([0,T],L^\infty(\R^d))$, we also get that
\begin{equation}\label{eq:T322}
\mathcal{T}_{3,2,2}(t)\,\leq\,\textcolor{black}{\mathcal{H}_\ve(t)}\,+\,o_{\ve\rightarrow 0}(1).
\end{equation}
\textcolor{black}{We precise that these two last estimates \eqref{eq:T321} and \eqref{eq:T322} are not uniform in $T$ in general since they involve norms in $L^\infty([0,T],H^2(\R^d))$ of the macroscopic quantities, which may not be uniform.}
\paragraph{\textcolor{black}{Modulated energy} estimate. } Finally, putting together the estimates \eqref{eq:T1}--\eqref{eq:T322}, we get that there exists a positive constant $C_T>0$ such that for all $t\in[0,T]$,
$$
\textcolor{black}{\mathcal{H}_\ve(t)}\,\leq\,\textcolor{black}{\mathcal{H}_\ve(0)} \,+\, C_T\left( \dfrac{1}{\ve^4}\|\rho_0^\ve-\rho_0\|_{L^2(\R^d)}^2 \,+\, \ve^{2/(d+6)} \,+\, o_{\ve\rightarrow 0}(1) \,+\, \ds\int_0^t\textcolor{black}{\mathcal{H}_\ve(s)} \,\md s \right).
$$
According to the assumptions \eqref{eq:ass3} and \eqref{eq:ass4} satisfied by the initial conditions, we get
$$\textcolor{black}{\mathcal{H}_\ve(t)}\,\leq\,o_{\ve\rightarrow 0}(1) \,+\, C_T\ds\int_0^t\textcolor{black}{\mathcal{H}_\ve(s)} \,\md s.
$$
Therefore, Grönwall's inequality yields that
\begin{equation}\label{eq:relentr}
\underset{\ve\rightarrow 0}{\lim}\,\underset{t\in[0,T]}{\sup}\,\textcolor{black}{\mathcal{H}_\ve(t)} \,=\, 0.
\end{equation}

\subsection{Conclusion}

Finally, let us conclude the proof of Theorem \ref{thm:mainresult} using the \textcolor{black}{modulated energy} estimate \eqref{eq:relentr} established in the previous subsection. Let $T>0$. We want to prove that the weak solution of the \textcolor{black}{transport} equation \eqref{eq:f} converges towards a monokinetic distribution as $\ve$ vanishes. First, we set
$$F^\ve(t,\mx,w)\,:=\,\ds\int f^\ve(t,\mx,v,w)\,\md v, \quad F^\ve(0,\mx,v)\,=\, F_0(\mx,w)\,:=\,\ds\int f_0^\ve(\mx,v,w)\,\md v.$$
Let us notice that since $f^\ve$ is compactly supported in $v$ for any $\ve>0$, we can choose a test function in \eqref{eq:weak} independent of $v\in\R$, so that the ditribution $F^\ve$ satisfies the following equation for all $\varphi\in\mathscr{C}^\infty_c([0,T)\times\R^{d+1})$:
$$\ds\int_0^T\int_{\R^{d+1}} \left( F^\ve\,\partial_t\varphi + \tau\,\left[ \ds\int_\R v\,f^\ve\,\md v - \gamma\,w\,F^\ve \right]\,\partial_w\varphi  \right)\,\md\mx\,\md w\,\md t \,+\,\ds\int_{\R^{d+1}}F_0^\ve\,\varphi(0) \,\md\mx\,\md w \,=\,0  ,$$
which is equivalent to satisfying for all $\varphi\in\mathscr{C}^1_c([0,T)\times\R^{d+1})$ the equation
\begin{multline}\label{eq:Feps}
\ds\int_0^T\int_{\R^{d+1}} F^\ve\,\left[ \partial_t\varphi + A\left( V(t,\mx),w \right)\,\partial_w\varphi \right]\,\md\mx\,\md w\,\md t \,+\,\ds\int_{\R^{d+1}}F_0^\ve\,\varphi(0) \,\md\mx\,\md w    \\
=\,\tau\ds\int_0^T\int_{\R^{d+2}}(V(t,\mx)-v)\,f^\ve\,\partial_w\varphi\,\md v \,\md w\,\md \mx\,\md t  ,
\end{multline}
where $V$ is solution to the second equation in \eqref{eq:lim2}. On the one hand, since $(F^\ve)_{\ve>0}$ in uniformly bounded by $1$ in $L^\infty([0,T],L^1(\R^{d+1}))$, we get that it converges weakly-$\ast$ \textcolor{black}{up to extraction} in $\mathcal{M}((0,T)\times\R^{d+1})$ towards a limit $F\in\mathcal{M}((0,T)\times\R^{d+1})$. Thus, we can pass to the limit on the left hand side of \eqref{eq:Feps} by linearity. On the other hand, from the \textcolor{black}{dissipation} estimate in \textcolor{black}{Proposition \ref{prop:dissip}} and the \textcolor{black}{modulated energy} estimate \eqref{eq:relentr}, we get that
\begin{eqnarray}\label{eq:int}
\ds\int_0^T\int f^\ve\,|v-V(t,\mx)|^2\,\md\mx\,\md v\,\md w\,\md t &\leq& 2\ds\int_0^T\int f^\ve\,\left( |v-V^\ve(t,\mx)|^2+|V^\ve(t,\mx)-V(t,\mx)|^2 \right)\,\md\mx\,\md v\,\md w\,\md t   \nonumber \\
&\longrightarrow& 0,  
\end{eqnarray}
as $\ve\rightarrow 0$. Consequently, since $\|\rho_0^\ve\|_{L^1}=1$, it yields with Cauchy-Schwarz inequality that:
$$\left|\ds\int_0^T\int_{\R^{d+2}}(V(t,\mx)-v)\,f^\ve\,\partial_w\varphi\,\md v \md w\,\md \mx\,\md t \right| \,\leq\, T^{1/2}\,\|\partial_w\varphi\|_\infty\,\ds\int_0^T\int|V(t,\mx)-v|^2\,f^\ve\,\md v \md w\,\md \mx\,\md t \,\longrightarrow\, 0,$$
as $\ve\rightarrow0$. Therefore, passing to the limit $\ve\rightarrow0$ in \eqref{eq:Feps}, it proves that $(V,F)$ is a solution of the system \eqref{eq:lim}. Furthermore, by uniqueness of the solution of \eqref{eq:lim}, we get the convergence of the sequence $(F^\ve)_{\ve>0}$.

Now, let us prove that for any $\varphi\in\mathscr{C}^0_b(\R^{d+2})$,
$$\ds\int\varphi(\mx,v,w)\,f^\ve(t,\mx,v,w)\,\md\mx\,\md v\,\md w\,\longrightarrow\,\ds\int\varphi(\mx,V(t,\mx),w)\,F(t,\md\mx, \md w),$$
strongly in $L^1_{\text{loc}}(0,T)$ as $\ve\rightarrow 0$. We start with proving that for all $0<t<t'\leq T$, for all $\varphi\in\mathscr{C}^1_c(\R^{d+2})$, 
\begin{equation}
\ds\int_t^{t'}\int f^\ve(s,\mx,v,w)\,\varphi(\mx,v,w)\,\md v\,\md w\,\md\mx\,\md s  \,\underset{\ve\rightarrow 0}{\longrightarrow}\, \ds\int_t^{t'}\int \varphi(\mx,V(s,\mx),w)\,F(s,\md\mx,\md w)\,\md s,
\end{equation} 
where $(V,F)$ is the solution on $[0,T]$ of the reaction-diffusion system \eqref{eq:lim} provided by Proposition \ref{prop:eqlim}, and we conclude using a density argument. Let $0<t<t'\leq T$. We can compute:
$$\mathcal{I}\,:=\,\left|\ds\int_t^{t'}\left(\int_{\R^{d+2}} f^\ve(s,\mx,v,w)\,\varphi(\mx,v,w)\,\md v\,\md w\,\md\mx \,-\, \int_{\R^{d+1}}\varphi(\mx,V(s,\mx),w) \,F(s,\md\mx,\md w)\right)\md s\right| \,\leq\, \mathcal{I}_1 \,+\, \mathcal{I}_2,$$
where 
$$\left\{\begin{array}{l}
\mathcal{I}_1 \,:=\, \ds\int_t^{t'}\int f^\ve(s,\mx,v,w)\,\left| \varphi(\mx,v,w)-\varphi(\mx,V(s,\mx),w) \right|\,\md v\,\md w\,\md\mx\,\md s, \\ \, \\
\mathcal{I}_2 \,:=\, \ds\int_t^{t'}\int \left| \varphi(\mx,V(s,\mx),w) \right|\,\left|F^\ve(s,\md\mx,\md w) - F(s,\md\mx,\md w)\right|\,\md s.
\end{array}\right.$$
On the one hand, using Cauchy-Schwarz inequality, we have
\begin{align*}
\mathcal{I}_1 \,&\leq\, \|\partial_v\varphi\|_\infty\,\ds\int_t^{t'}\int f^\ve(s,\mx,v,w)\,\left| v-V(s,\mx) \right|\,\md v\,\md w\,\md\mx\,\md s    \\
&\leq\,\|\partial_v\varphi\|_\infty\,\left( \ds\int_0^T\int f^\ve(s,\mx,v,w)\,\md v\,\md w\,\md\mx\,\md s \right)^{1/2}\,\left( \ds\int_0^T\int f^\ve(s,\mx,v,w)\,|v-V(s,\mx)|^2\,\md v\,\md w\,\md\mx\,\md s \right)^{1/2}  \\
&=\, \|\partial_v\varphi\|_\infty\,T^{1/2}\,\left( \ds\int_0^T\int f^\ve(s,\mx,v,w)\,|v-V(s,\mx)|^2\,\md v\,\md w\,\md\mx\,\md s \right)^{1/2}  .
\end{align*}
Consequently, using the convergence from \eqref{eq:int}, we get
\begin{equation}\label{eq:I1}
\underset{\ve\rightarrow 0}{\lim}\,\mathcal{I}_1 \,=\, 0.
\end{equation}
On the other hand, the second term $\mathcal{I}_2$ converges to zero as $\ve$ goes to zero since $(F^\ve)_{\ve>0}$ converges weakly-$\ast$ towards $F$ in $\mathcal{M}((0,T)\times\R^{d+1})$. Consequently, we can conclude that
\begin{equation}
\underset{\ve\rightarrow 0}{\lim}\,\mathcal{I} \,=\, 0.
\end{equation}
Using a density argument, this shows the convergence of $f^\ve$ in $L^1_\text{loc} ((0,T),\mathcal{M}(\R^{d+2}))$ towards a monokinetic distribution, which concludes the proof of Theorem \ref{thm:mainresult}.

\section{Proof of Proposition \ref{prop:eqlim}}\label{sec:eqlim}

This subsection is devoted to the proofs of Proposition \ref{prop:eqlim} and its Corollary \ref{cor:eqlim}, that is to the construction of a solution to the system \eqref{eq:lim2}. The main difficulty lies in the fact that the function $\rho_0$ can reach $0$, so the second equation in \eqref{eq:lim2} is not well-defined on $\R^d$. A solution to overcome this problem is to construct a solution to \eqref{eq:lim2} from a weak solution of the reaction-diffusion \textcolor{black}{FHN} system \eqref{eq:FHNlocal} in the sense of Definition \ref{defi:weaksol}.

\textcolor{black}{Furthermore, we have to be especially careful} to prove the existence and uniqueness of a weak solution to the reaction-diffusion system \eqref{eq:FHNlocal} \textcolor{black}{since} it is not a parabolic system. A way to circumvent this issue is to consider for all $\delta\geq 0$ the linear operator:
\begin{equation}
\mathcal{L}_{\rho_0+\delta}:V \,\mapsto\, \sigma \big[ (\rho_0+\delta)\,\Delta_\mx V  \,+\, 2\,\nabla_\mx \rho_0\,\cdot\,\nabla_\mx V \big].
\end{equation}

First of all, for any given initial data $V_0\in L^2(\R^d)$ and for all $\de>0$, we prove the existence and uniqueness of $V_\de$ a weak solution of the approximated parabolic system
\begin{equation}\label{eq:FHNlocal2}
\partial_t V_\de \,+\,\mathcal{L}_{\rho_0+\de}(V_\de) \,=\, N(V_\de) \,-\, W[V_\de],
\end{equation}
where for all $V:\R^+\times\R^d\rightarrow\R$, we define 
\begin{equation}\label{eq:W[V]}
W[V]:(t,\mx) \,\mapsto\,  e^{-\tau\,\gamma\,t}W_0(\mx) \,+\,\tau\ds\int_0^t e^{-\tau\,\gamma\,(t-s)}V(s,\mx)\,\md s,
\end{equation}
where $W_0:\R^d\rightarrow\R$ is a given initial data in $H^2(\R^d)$. 
\textcolor{black}{We will be able to pass to the limit $\de\rightarrow0$ thanks to \textit{a priori} estimates of the $H^2$ norm of the solution of \eqref{eq:FHNlocal2} which are uniform in $\de$.}

In the rest of this article, for all $k\in\{0,1,2\}$, we note $\langle\cdot,\cdot\rangle_{H^k(\R^d)}$ the scalar product of $H^k(\R^d)$ defined as follows:
$$\langle U,V\rangle_{H^k(\R^d)} \,:=\, \underset{\alpha\in\N^d,|\alpha|\leq k}{\ds\sum} \, \ds\int \partial^\alpha U \, \partial^\alpha V \, \md\mx,$$
for all $U,V\in H^k(\R^d)$, where for all $\alpha=(\alpha_1,...,\alpha_d)\in\N^d$, $\partial^\alpha \,=\, \partial_{\mx_1}^{\alpha_1}\,...\,\partial_{\mx_d}^{\alpha_d}\,$.

\subsection{\textit{A priori} estimates}\label{subsec:apriori}

\textcolor{black}{The purpose of this subsection is to derive an \textit{a priori} estimate of the $H^2$ norm of a weak solution to the \textcolor{black}{FHN} reaction-diffusion system \eqref{eq:FHNlocal2} uniform in $\de$.}
\begin{lem}\label{lem:apriori}
Let $T>0$. Consider an initial data $\rho_0$ satisfying \eqref{hyp:rho0} and $V_0\in H^2(\R^d)$ satisfying \eqref{hyp:V0}. Let $\de\geq 0$. Assume that there exists 
$$V_\de\in L^\infty([0,T],H^2(\R^d))\,\cap\,\mathscr{C}^0([0,T],L^2(\R^d))$$
a weak solution of the reaction-diffusion equation for $t>0$ and $\mx\in\R^d$:
\begin{equation}\label{eq:approx1}
\partial_t V_\de \,-\,  \mathcal{L}_{\rho_0+\delta} V_\de  \,=\, S,
\end{equation}
where $S\in L^\infty([0,T],H^2(\R^d))$ are two general source terms. Then, there exists a positive constant $C>0$ independent of $\de$ such that for all $t\in [0,T]$
\begin{equation}\label{eq:Hk}
\|V_\de(t)\|_{H^{2}(\R^d)}^2 \,+\,\de\ds\int_0^t\|V_\de(s)\|_{H^{3}(\R^d)}^2    \,\leq\, \|V_0\|^2_{H^{2}(\R^d)} \,+\, C\ds\int_0^t\left(\|V_\de(s)\|_{H^{2}(\R^d)}^2 \,+\, \langle S(s) \,,\,V_\de(s)\rangle_{H^{2}(\R^d)} \right)\,\md s .
\end{equation}
\end{lem}

\begin{Proof}
We postpone the proof to Appendix \ref{app:lemapriori}.
\end{Proof}


\begin{cor}\label{cor:apriori_N(V)}
Let $\de\geq 0$. Consider an initial data $\rho_0$ satisfying \eqref{hyp:rho0} and $(V_0,W_0)$ satisfying \eqref{hyp:V0}. Let $T>0$ such that there exists 
$$V_\de \in  L^\infty([0,T],H^2(\R^d))\,\cap\,\mathscr{C}^0([0,T],L^2(\R^d))$$
a weak solution of the reaction-diffusion system \eqref{eq:FHNlocal2}. Then, there exists a finite constant $C_{T}>0$ independant of $\de$ such that for all $t\in[0,T]$,
\begin{equation}
\|V_\de(t)\|_{H^2(\R^d)} \,\leq \,C_{T}.
\end{equation}
\end{cor}

\begin{Proof}
We postpone the proof to Appendix \ref{app:lemapriori_N(V)}.
\end{Proof}

\subsection{Case of a positive $\de$}

Let $\de>0$. \textcolor{black}{This subsection focuses on the existence and uniqueness of the reaction-diffusion equation \eqref{eq:FHNlocal2}.}

\begin{lem}\label{lem:eqdelta}
Consider an initial data $\rho_0$ satisfying \eqref{hyp:rho0} and $V_0\in H^2(\R^d)$. Then, for all $T>0$ and for all $\de>0$, there exists a unique weak solution $V_\de$ of the diffusion equation
\begin{equation}\label{eq:approx10}
\partial_t V_\de \,-\, \mathcal{L}_{\rho_0+\delta} V_\de \,=\, N(V_\de) \,-\, W[V_\de],
\end{equation}
such that
$$V_\de\in L^\infty([0,T],H^2(\R^d))\cap L^2([0,T],H^3(\R^d))\cap\mathscr{C}^0([0,T],L^2(\R^d)),$$
and $V_\de$ satisfies the energy estimate \eqref{eq:Hk} with $S=N(V_\de)-W[V_\de]$.
\end{lem}

\begin{Proof}
The proof relies on classical methods explained in \textcolor{black}{Section 7.1 in} \cite{EVA}. According to Lemma \ref{lem:apriori}, $V_\de$ satisfies the energy estimate \eqref{eq:Hk} with $S=N(V_\de)-W[V_\de]$. 
\end{Proof}

\subsection{Proof of Proposition \ref{prop:eqlim}}

{\color{black}\paragraph{Step 1: Existence.}Now}, let us pass to the limit $\de\rightarrow0$ in the approximated equation \eqref{eq:FHNlocal2}, to prove Proposition \ref{prop:eqlim}. Let $V_0$ and $W_0\in H^2(\R^d)$. For all $\de>0$, Lemma \ref{lem:eqdelta} yields the existence of a weak solution $V_\de$ to the reaction-diffusion equation \eqref{eq:FHNlocal2}. According to Corollary \ref{cor:apriori_N(V)}, there exists a positive constant $K_1>0$ independent of $\de$ such that for all $\de>0$,
$$\|V_\de\|_{L^\infty([0,T],H^2(\R^d))}\,\leq\, K_1.$$
Let $(\de_n)_{n\in\mathbb{N}}$ be a sequence of positive reals such that $\delta_n\rightarrow 0$ as $n\rightarrow\infty$. Therefore, there exists a function $V\in L^\infty([0,T],H^2(\R^d))$ such that up to extraction,
$$V_{\de_n} \,\rightharpoonup\, V, $$
weakly-$\star$ in $L^\infty([0,T],H^2(\R^d))$ as $\de\rightarrow 0$. \textcolor{black}{Let us prove that the sequence $(V_{\de_n})_{n\in\mathbb{N}}$ also converges in a strong sense using Arzelà-Ascoli Theorem. On the one hand, since for all $n\in\mathbb{N}$ the function $V_{\de_n}$ is a weak solution of \eqref{eq:approx10}, then there exists a constant $K_2>0$ independent of $n$ such that for all $n\in\mathbb{N}$,
$$\left\|\partial_t V_{\de_n}\right\|_{L^\infty([0,T],L^2(\R^d))} \,\leq\,K_2.$$
On the other hand, for all $t\in[0,T]$, the set $\{V_{\de_n}(t) \,|\, n\in\mathbb{N}\}$ is not relatively compact in $L^2(\R^d)$. In order to validate this last assumption we have to restrict the domain of integration to open bounded sets. Indeed, for all open bounded subset $\Omega$ of $\R^d$, the set $\{V_{\de_n}(t)|_{\Omega} \,|\, n\in\mathbb{N}\}$ is relatively compact in $L^2(\Omega)$ since the inclusion $H^2(\Omega)\,\subset\,L^2(\Omega)$ is compact. In the following, for all $r>0$, we note $B_r$ the ball of $\R^d$ of radius $r$ centered in $0$.
}

\textcolor{black}{For all $k\in\mathbb{N}^*$ and all $n\in\mathbb{N}$, let us consider $V_{n,k}\,:=\,V_{\de_n}|_{B_k}$. We use a diagonal extraction argument and the Arzelà-Ascoli theorem to obtain that for all $k\in\mathbb{N}^*$, there exists an extraction $\phi_k$ such that for all $k'\leq k$, the sequence $(V_{\phi_{k}(n),k'})_{n\in\mathbb{N}}$ converges towards a function $V_{\infty,k'}$ strongly in $\mathscr{C}^0\left(([0,T],L^2(B_{k'})\right)$ as $n$ goes to infinity.}

\textcolor{black}{We claim that for all $k\in\mathbb{N}^*$, $V_{\infty,k}\,=\, V_{\infty,k+1}|_{B_k}$. Indeed, we have
\begin{align*}
&\|V_{\infty,k}\,-\, V_{\infty,k+1}\|_{L^\infty\left([0,T],L^2(B_k)\right)} \\
&\leq\, \|V_{\infty,k}\,-\, V_{\phi_{k+1}(n),k}\|_{L^\infty\left([0,T],L^2(B_k)\right)} \,+\, \|V_{\infty,k+1}\,-\, V_{\phi_{k+1}(n),k+1}\|_{L^\infty\left([0,T],L^2(B_k)\right)} \\ 
&\quad \,+\, \|V_{\phi_{k+1}(n),k}\,-\, V_{\phi_{k+1}(n),k+1}\|_{L^\infty\left([0,T],L^2(B_k)\right)}.
\end{align*}
Furthermore, we also have that $V_{\phi_{k+1}(n),k}\,=\, V_{\phi_{k+1}(n),k+1}|_{B_k}$ according to the definition of the sequence $(V_{n,k})_{n\in\mathbb{N}}$. Thus, passing to the limit $n\rightarrow+\infty$ in the last inequality, we can prove our claim.}

\textcolor{black}{Consequently, there exists a function $V_\infty\in\mathscr{C}([0,T],L^2_\text{loc}(\R^d))$ such that for all $k\in\mathbb{N}^*$, the function $V_{\infty,k}$ is the restriction of $V_k$ to the domain $B_k$. Thus, the sequence $(V_{\phi_n(n),n})_{n\in\mathbb{N}}$ strongly converges towards $V_\infty$ in $L^\infty([0,T],L^2_\text{loc}(\R^d))$. Therefore, $V=V_\infty$, and thus,
$$V\in L^\infty\left([0,T],H^2(\R^d)\right)\,\cap\,\mathscr{C}^0\left([0,T],L^2_\text{loc}(\R^d)\right).$$
Since the sequence $(V_{\phi_n(n),n})_{n\in\mathbb{N}}$ strongly converges towards $V$, we can pass to the limit in the weak formulation of the equation \eqref{eq:FHNlocal2} when the space of test functions is $\mathscr{C}^\infty_c(\R^{d+2})$. Consequently, $V$ is a solution in the sense of distributions of \eqref{eq:FHNlocal} with $\de=0$. Furthermore, since $V\in L^\infty([0,T],H^2(\R^d))$, we can deduce that $\partial_t V\in L^\infty([0,T],L^2(\R^d))$, and thus
$$V\in W^{1,\infty}([0,T],L^2(\R^d)).$$
Therefore, using classical arguments detailed in the paragraph 5.9.2 from \cite{EVA}, we get:
$$V\in\mathscr{C}^0([0,T],H^1(\R^d)).$$
Furthermore, since the space of test functions $\mathscr{C}^\infty_c(\R^{d+2})$ is dense in $H^1(\R^{d+2})$, we get that $V$ is a weak solution of \eqref{eq:FHNlocal} in the sense of Definition \ref{defi:weaksol}.
}

{\color{black}
{\color{black}\paragraph{Step 2: Uniqueness.}Finally}, let us justify the uniqueness of the solution. Let $V$ and $\widetilde{V}$ be two solutions of the Cauchy problem \eqref{eq:FHNlocal}. If we consider the function $V-\widetilde{V}$, we notice that it satisfies the equation \eqref{eq:approx1} with $\de=0$ and with the source term
$$S \,=\, N(V)-N(\widetilde{V}) \,-\, \left( W[V] - W[\widetilde{V}] \right).$$
Let us estimate the $L^2$ norm of $V-\widetilde{V}$. Following the same computations as in Lemma \ref{lem:apriori}, we find that there exists a positive constant $C>0$ such that for all $t\in[0,T]$,
\begin{equation}\label{eq:L2}
\|(V-\widetilde{V})(t)\|_{L^2(\R^d)} \,\leq\, C\ds\int_0^t\left( \|(V-\widetilde{V})(s)\|_{L^2(\R^d)} \,+\, \langle S(s) \, , \, (V-\widetilde{V})(s) \rangle_{L^2(\R^d)} \right)\,\md s .
\end{equation}
It remains to estimate the scalar product. According to the assumption \eqref{hyp:N} satisfied by the nonlinearity $N$, we get that there exists a positive constant $C_T>0$ such that for all $s\in[0,t]$,
$$\left\{\begin{array}{l}
\ds\int (V-\widetilde{V})(s,\mx) \, (N(V)-N(\widetilde{V}))(s,\mx) \,\md \mx \,\leq\,\kappa_2\,\|(V-\widetilde{V})(s)\|_{L^2(\R^d)}^2, \\ \, \\
\ds\int (V-\widetilde{V})(s,\mx) \, (W[V]-W[\widetilde{V}])(s,\mx) \,\md \mx \,\leq\,C_T\,\|V-\widetilde{V}\|_{L^\infty([0,s],L^2(\R^d))}^2.
\end{array}\right.$$
Then, taking the supremum over time in \eqref{eq:L2}, we get that there exists a positive constant $C_T>0$ such that for all $t\in[0,T]$:
$$\|V-\widetilde{V}\|_{L^\infty([0,t],L^2(\R^d))} \,\leq\, C_T\ds\int_0^t\|V-\widetilde{V}\|_{L^\infty([0,s],L^2(\R^d))} \,\md s .$$
According to Grönwall's inequality, we get that $V=\widetilde{V}$.
}

\subsection{Conclusion: proof of Corollary \ref{cor:eqlim}}


Let $T>0$. Let us denote with $(\widetilde{V},\widetilde{W})$ the weak solution of the equation \eqref{eq:FHNlocal} provided by Proposition \ref{prop:eqlim} with initial condition $(\rho_0,V_0,W_0)$. Our proof is organised in two steps. First of all, we claim that for all solution $(V,F)$ of the system \eqref{eq:lim}, the two functions $V$ and $\widetilde{V}$ coincide almost everywhere on $[0,T]\times\R^d$. Then, we prove the existence of a measure solution $F$ so that $(\widetilde{V},F)$ is a solution of the system \eqref{eq:lim}. 

Before starting the proof, notice that for all $\mx\in\R^d$ such that $\rho_0(\mx)=0$, the system \eqref{eq:FHNlocal} reduces to the following system of ODEs
\begin{equation}\label{eq:ODEs}
\left\{\begin{array}{l}
\partial_t \widetilde{V} \,=\, N(\widetilde{V}) - \widetilde{W}, \\ \, \\
\partial_t \widetilde{W} \,=\, A(\widetilde{V},\widetilde{W}).
\end{array}\right.
\end{equation}
\textcolor{black}{For almost every $\mx\in\R^d$ such that $\rho_0(\mx)=0$, since $V_0(\mx)\,=\,W_0(\mx)\,=\,0$ according to \eqref{def:rho0W0}--\eqref{hyp:V0=0}, one can directly conclude that for all $t\in[0,T]$,}
$$\widetilde{V}(t,\mx)\,=\,\widetilde{W}(t,\mx)\,=\,0.$$

\paragraph{Step 1: Uniqueness.} Now, let us prove that for any solution $(V,F)$ on $[0,T]$ of \eqref{eq:lim} in the sense of Definition \ref{defi:sol}, the function $V$ coincides with $\widetilde{V}$ on $[0,T]\times\R^d$. Suppose that $(V,F)$ is a solution of \eqref{eq:lim} such that $F$ has a finite second moment in $w$. If we define for all $(t,\mx)\in[0,T]\times\R^d$
$$\rho_0(\mx)\,W(t,\mx)\,:=\,\ds\int w\,F(t,\mx,\md w),$$
then the \textcolor{black}{triple} $(\rho_0,\rho_0\,V,\rho_0\,W)$ satisfies the reaction-diffusion equation \eqref{eq:rhoFHNlocal}. 
On the one hand, by definition, for almost every $\mx\in\R^d$ such that $\rho_0(\mx)=0$, for all $t\in[0,T]$, 
$$V(t,\mx)\,=\,\widetilde{V}(t,\mx)\,=\,0\, , \quad W(t,\mx)\,=\,\widetilde{W}(t,\mx)\,=\,0.$$
According to the notion of solution of \eqref{eq:lim} from Definition \ref{defi:sol}, $V\in L^\infty([0,T],H^2(\R^d))$, so $\Delta_\mx V(t,\mx)$ is defined for almost every $\mx\in\R^d$. Consequently, the \textcolor{black}{pair} $(V,W)$ satisfies \eqref{eq:FHNlocal} pointwise for all $t\in[0,T]$ and almost every $\mx\in\R^d$ such that $\rho_0(\mx)=0$.

On the other hand, for all $\mx\in\R^d$ such that $\rho_0(\mx)>0$, the equation \eqref{eq:rhoFHNlocal} reduces to the reaction-diffusion system \eqref{eq:FHNlocal}. Therefore, $(V,W)$ satisfies the equation \eqref{eq:FHNlocal} for all $t\in[0,T]$ and almost every $\mx\in\R^d$, with initial condition $(V_0,W_0)$. Consequently, the \textcolor{black}{pair} $(V,W)$ satisfies the reaction-diffusion system \eqref{eq:FHNlocal} in the sense of Definition \ref{defi:weaksol}. Thus, by uniqueness of the solution of the equation \eqref{eq:FHNlocal}, we can conclude that 
$$(V,W)\,=\,(\widetilde{V},\widetilde{W}).$$

\paragraph{Step 2: Existence. }Then, let us prove the existence of a measure solution of the first equation in \eqref{eq:lim}. With $V$ the first component of the solution of the system \eqref{eq:FHNlocal}, let us consider the transport equation for $t>0$, $\mx\in\R^d$ and $w\in\R$:
\begin{equation}\label{eq:F}
\left\{\begin{array}{l}
\partial_t F(t,\mx,w) \,+\,\partial_w\left( A(V(t,\mx),w)\,F(t,\mx,w) \right) \,=\,0,   \\ \, \\
F|_{t=0}\,=\,F_0.
\end{array}\right.
\end{equation}
In order to solve \eqref{eq:F}, we introduce the associated system of characteristic curves for all $(s,t,\mx,w)\in[0,T]^2\times\R^{d+1}$:
\begin{equation}\label{eq:W}
\left\{\begin{array}{l l}
\dfrac{\md}{\md s}\mathcal{W}(s) \,=\, A(V(s,\mx),\mathcal{W}(s)),&   \\ \,  \\
\mathcal{W}(t)\,=\,w.
\end{array}\right.
\end{equation}
Since the function $A$ grows linearly with respect to $w$, and the function $V$ is regular enough, we get the global existence and uniqueness of a solution of the characteristic equation \eqref{eq:W}{\color{black}. Then, using the theory of characteristics, we get the existence of a unique solution to the transport equation \eqref{eq:F}. Then, we directly get from \eqref{eq:F} and the assumptions \eqref{mom2} and \eqref{def:rho0W0} that there exists a positive constant $C_T>0$ such that for all $t\in[0,T]$ and all $\mx\in\R^d$,
$$\int\!\!\int_{\R^{d+1}} |w|^2\,F(t,\md\mx,\md w)\,\leq\,C_T, \quad \dfrac{1}{\rho_0(\mx)}\ds\int_\R w\, F(t,\mx,\md w) \,=\, W(t,\mx).$$
}
Consequently, the unique solution of the system \eqref{eq:lim} is $(V,F)$ where $(V,W)$ is the weak solution of \eqref{eq:FHNlocal} and $F$ is the unique solution of the transport equation \eqref{eq:F}.
\textcolor{black}{\section{Acknowledgements}
JC acknowledges support from an ANITI (Artificial and Natural Intelligence Toulouse Institute) Research Chair.}

\appendix

\section{Proof of Lemma \ref{lem:apriori}}\label{app:lemapriori}
Our approach consists in studying the variations of $\|V_\de\|_{L^2(\R^d)}$ and $\|\Delta_\mx V_\de\|_{L^2(\R^d)}$, in order to conclude with interpolations. Let $t\in[0,T]$. For the sake of simplicity, in the rest of this proof, we note $V$ instead of $V_\de$. First of all, we get:
\begin{align*}
\dfrac{1}{2}\,\dfrac{\md}{\md t}\|V\|_{L^2(\R^d)}^2   \,&=\, \ds -\sigma \left( \int  (\rho_0+\de)\,\left|\nabla_\mx V\right|^2 \,\md\mx \,-\,\int  (\nabla_\mx \rho_0\,\cdot\,\nabla_\mx V) \,V\,\md\mx  \right) \,+\,\int  S\,V\,\md\mx  \\
&= \, \ds - \sigma \left( \int (\rho_0+\de)\,\left|\nabla_\mx V\right|^2\,\md\mx \,+\,\dfrac{1}{2}\int  \Delta_\mx \rho_0\,|V|^2\,\md\mx  \right) \,+\,\int  S\,V\,\md\mx .
\end{align*}
Thus, \textcolor{black}{since $\rho_0\in\mathscr{C}^3_b(\R^d)$}, we can conclude that there exists a constant $C>0$ such that
\begin{equation}\label{estimL2}
\dfrac{1}{2}\,\dfrac{\md}{\md t}\|V\|_{L^2(\R^d)}^2  \,+\, \sigma \ds\int  (\rho_0+\de)\,\left|\nabla_\mx V\right|^2\,\md\mx \,\leq\, C\,\|V\|_{L^2(\R^d)}^2 \,+\,\int  S\,V\,\md\mx .
\end{equation}
Then, we also get that
\begin{align*}
&\dfrac{1}{2}\,\dfrac{\md}{\md t} \|\Delta_\mx V\|_{L^2(\R^d)}^2  \\
&=\,-\ds\int\Delta_\mx((\rho_0+\de)\,\nabla_\mx V)\cdot\nabla_\mx
\Delta_\mx V\,\md\mx   \,+\, \ds\int\Delta_\mx(\nabla_\mx\rho_0\cdot\nabla_\mx V)\,\Delta_\mx V\,\md\mx \,+\,\ds\int \Delta_\mx S\,\Delta_\mx V\,\md\mx   \\
&=\,-\ds\int\left[\Delta_\mx \rho_0\,\nabla_\mx V + 2\,\Delta_\mx V \nabla_\mx\rho_0 + (\rho_0+\de)\,\nabla_\mx\Delta_\mx V \right]\cdot\nabla_\mx\Delta_\mx V \,\md\mx   \\
&~~~~+\ds\int\left[ \nabla_\mx\Delta_\mx \rho_0\cdot\nabla_\mx V + 2\, \Delta_\mx\rho_0\,\Delta_\mx V + \nabla_\mx\rho_0\cdot\nabla_\mx\Delta_\mx V \right]\,\Delta_\mx V\,\md\mx   \,+\, \ds\int \Delta_\mx S\,\Delta_\mx V \,\md\mx.
\end{align*}
Using Green's formula on the term $\ds\int \Delta_\mx\rho_0\,\nabla_\mx V\cdot\nabla_\mx\Delta_\mx V\,\md\mx$, we compute that
\begin{align*}
\dfrac{1}{2}&\,\dfrac{\md}{\md t}\|\Delta_\mx V\|_{L^2(\R^d)}^2  \\
&=\,-\ds\int(\rho_0+\de)\,|\nabla_\mx\Delta_\mx V|^2\,\md\mx \,+\,3\int\Delta_\mx\rho_0\,\left|\Delta_\mx V\right|^2\,\md\mx \,+\, \ds\int\nabla_\mx\Delta_\mx\rho_0\cdot\nabla_\mx V\,\Delta_\mx V\,\md\mx   \\
&~~~~+\ds\int \Delta_\mx S\,\Delta_\mx V \,\md\mx\\
&\leq\,-\ds\int(\rho_0+\de)\,|\nabla_\mx\Delta_\mx V|^2\,\md\mx \,+\, 3\,\|\rho_0\|_{\mathscr{C}^2(\R^d)}\,\|\Delta_\mx V\|_{L^2(\R^d)}^2 \,+\, \|\rho_0\|_{\mathscr{C}^3(\R^d)}\,\|\nabla_\mx V\|_{L^2(\R^d)}\,\|\Delta_\mx V\|_{L^2(\R^d)}   \\
&~~~~+\ds\int \Delta_\mx S\,\Delta_\mx V \,\md\mx\\
&\leq\,-\ds\int(\rho_0+\de)\,|\nabla_\mx\Delta_\mx V|^2\,\md\mx \,+\, C\left( \|V\|_{L^2(\R^d)}^2 \,+\, \|\Delta_\mx V\|_{L^2(\R^d)}^2 \right) \,+\ds\int \Delta_\mx S\,\Delta_\mx V \,\md\mx,
\end{align*}
where $C>0$ is a positive constant. Consequently, we get that for all $t\in[0,T]$,
\begin{multline}\label{estimH2}
\dfrac{1}{2}\,\dfrac{\md}{\md t}\,\|\Delta_\mx V(t)\|^2_{L^{2}(\R^d)} \,+\, \textcolor{black}{\ds\int(\rho_0+\de)\,|\nabla_\mx\Delta_\mx V|^2\,\md\mx} \\
\leq\, C\,\left(\|\Delta_\mx V(t)\|^2_{L^2(\R^d)} \,+\,\|V(t)\|^2_{L^2(\R^d)}\right)  \,+\, \ds\int \Delta_\mx S(t)\,\Delta_\mx V(t)\,\md\mx.
\end{multline}
Finally, integrating the estimates \eqref{estimL2} and \eqref{estimH2} between $0$ and $t$ for $t\in[0,T]$, we get the estimate \eqref{eq:Hk} \textcolor{black}{since the $H^2$ norm is equivalent to the norm $\|\cdot\|_{L^2(\R^d)} + \|\Delta_\mx\cdot\|_{L^2(\R^d)}$.}

\section{Proof of Corollary \ref{cor:apriori_N(V)}}\label{app:lemapriori_N(V)}

For the sake of simplicity, in the rest of this section, we note $V$ instead of $V_\de$. According to Lemma \ref{lem:apriori}, the estimate \eqref{eq:Hk} holds with $S=N(V)-W[V]$. To obtain Corollary \ref{cor:apriori_N(V)} from the energy estimate \eqref{eq:Hk}, we need to estimate the scalar product $\langle V,N(V) \rangle_{H^2(\R^d)}$. First of all, for all $t\in[0,T]$, since $N$ satisfies the property \eqref{hyp:N}, we have:
$$\ds \int V \, N(V) \,\md\mx\,\leq\,\kappa_1\,\|V\|^2_{L^2(\R^d)} \,-\,\kappa_1'\|V\|^4_{L^4(\R^d)}\,\leq\, \kappa_1\,\|V\|^2_{L^2(\R^d)}.$$
\textcolor{black}{Then, we only give details of computation of the integral of the product of $\Delta_\mx V$ and $\Delta_\mx N(V)$. Using} Young's inequality for some small parameter $\theta >0$, we also obtain:
\begin{align*}
\ds\int \Delta_\mx V  \, \Delta_\mx N(V) \,\md\mx\,&=\,\ds\int  |\Delta_\mx V|^2\,\left[-3|V|^2 + 2(1+a)V - a\right]\,\md\mx \,+\,\int  \Delta_\mx V\,\left| \nabla_\mx V \right|^2\,\left[ -6 V + 2(1+a) \right]\,\md\mx   \\
&\leq\,-3\,\ds\int  |\Delta_\mx V|^2\,|V|^2\,\md\mx  \,-\,a\ds\int\left|\Delta_\mx V\right|^2\,\md\mx   \\
&~~~+\,\dfrac{(1+a)}{\theta}\ds\int|\Delta_\mx V|^2\,\md\mx \,+\, (1+a)\,\theta \int|\Delta_\mx V|^2\,|V|^2\,\md\mx   \\
&~~~ +\, \dfrac{3}{\theta}\,\ds\int\left|\nabla_\mx V\right|^4\,\md\mx  \,+\, 3\,\theta \ds\int |V|^2\,\left|\Delta_\mx V\right|^2\,\md\mx \\
&~~~+\, (1+a)\ds\int\left|\Delta_\mx V\right|^2\,\md\mx \,+\, (1+a)\ds\int\left|\nabla_\mx V\right|^4\,\md\mx.
\end{align*}
Consequently, if we consider $\theta$ small enough so that $(1+a)\,\theta\,+\,3\,\theta \,\leq\,3$, we obtain:
\begin{equation}\label{eq:app1}
\ds\int \Delta_\mx V  \, \Delta_\mx N(V) \,\md\mx\,\leq\, \left( 1+\dfrac{1}{\theta} \right)\,(1+a)\,\ds\int\left|\Delta_\mx V\right|^2\,\md\mx   \,+\,\left( 1+a+\dfrac{3}{\theta}  \right)\ds\int\left|\nabla_\mx V\right|^4\,\md\mx.
\end{equation}
Therefore, to conclude, we only need to find a uniform bound of the $L^4$ norm of $\nabla_\mx V(t)$. We apply the Gagliardo-Nirenberg inequality on $\|\nabla_\mx V_0\|_{L^4(\R^d)}$, which yields that there exists a positive constant $C>0$ such that:
$$\|\nabla_\mx V_0\|_{L^4(\R^d)} \,\leq\,C\,\|V_0\|_{H^2(\R^d)}^{\frac{d}{4}}\,\|\nabla_\mx V_0\|_{L^2(\R^d)}^{1-\frac{d}{4}}\,<\,+\infty,$$
Hence, since $V_0\in H^2(\R^d)$, we have $\nabla_\mx V_0\in L^4(\R^d)$, and it still holds if we replace $V_0$ with $W_0$. 
Furthermore, $\nabla_\mx V$ satisfies in the weak sense the following equation on $(0,T]\times\R^d$
$$ \partial_t\left( \nabla_\mx V \right) \,=\, \sigma \left[ \nabla_\mx\left( \rho_0\,\Delta_\mx V \right) \,+\,2\,\nabla_\mx\left( \nabla_\mx\rho_0\cdot\nabla_\mx V \right) \right] \,+\,\nabla_\mx V\,N'(V)\,-\,\nabla_\mx W[V] . $$
{\color{black}We can use this last equation and similar computations as before to estimate the $L^4$ norm of $\nabla_\mx V$. Thus, we get that there exists a positive constant $K_T>0$ such that
\begin{equation}\label{eq:app2}
\underset{t\in[0,T]}{\sup}\,\|\nabla_\mx V(t)\|_{L^4(\R^d)}^4 \,\leq\, K_T.
\end{equation}
}
Therefore, we conclude from \eqref{eq:app1}-\eqref{eq:app2} that there exists a positive constant $C$ such that \textcolor{black}{for all $t\in[0,T]$}:
\begin{align*}
\ds\int \Delta_\mx V  \, \Delta_\mx N(V) \,\md\mx\,&\leq\, C\,\|V\|_{H^2(\R^d)}^2 \,+\,K_T.
\end{align*}
\bibliography{plain}

\end{document}